\documentclass[draft,preprint,12pt,numbers,sort&compress]{elsarticle}
\usepackage{mathrsfs}
\usepackage{amssymb}
\usepackage{amsthm}
\usepackage{mathrsfs}
\usepackage[centertags]{amsmath}
\usepackage{amsfonts}

\usepackage{color}

\input amssym.def
\input amssym.tex

\date{}

\newtheorem{Theorem}{Theorem}[section]

\newtheorem{Lemma}{Lemma}[section]

\newcommand\R{\mbox{\bf R}}
\newcommand\N{\mbox{\bf N}}
\newcommand\SR{\mbox{\scriptsize\bf R}}

\newcommand{\definition}{{\lower .5ex
  \hbox{$\>\>\stackrel{\triangle}{=}\>\>$} }}
\newcommand\supp{\mathop{\rm supp}}
\newcommand\mes{\mathop{\rm mes}}

\begin{document}

\baselineskip=22pt
\thispagestyle{empty}

\bigskip

\begin{center}{\Large\bf The Cauchy problem for the Ostrovsky equation }\\[1ex]
{\Large\bf with negative dispersion at the critical regularity}\\[4ex]

{Yongsheng LI$^\dag$  \quad\quad Jianhua HUANG$^\ddag$\footnote{Corresponding Author: J.Huang,\quad Email: jhhuang32@nudt.edu.cn} \quad and \quad Wei YAN$^*$}\\[2ex]
{$^\dag$Department of Mathematics, South China University of Technology,}\\
{Guangzhou, Guangdong 510640, P. R. China}\\[2ex]
{$^\ddag$ College of Science, National University of Defense and Technology,}\\
{ Changsha, P. R. China\quad  410073}\\[2ex]
{$^*$College of Mathematics and Information Science, Henan Normal University,}\\
{Xinxiang, Henan 453007, P. R. China}
\end{center}

\noindent{\bf Abstract.}   In this paper, we investigate the Cauchy problem
for the Ostrovsky equation
\begin{eqnarray*}
      \partial_{x}\left(u_{t}-\beta \partial_{x}^{3}u
  +\frac{1}{2}\partial_{x}(u^{2})\right)
    -\gamma u=0,
\end{eqnarray*}
in the Sobolev space $H^{-3/4}(\R)$. Here $\beta>0(<0)$
corresponds to the positive (negative) dispersion of the media, respectively.
P. Isaza and  J. Mej\'{\i}a (J. Diff. Eqns. 230(2006), 601-681; Nonli. Anal. 70(2009), 2306-2316),
 K. Tsugawa (J. Diff. Eqns. 247(2009), 3163-3180) proved that the problem
 is locally well-posed in $H^s(\R)$ when
$s>-3/4$ and ill-posed when $s<-3/4$. By using some modified Bourgain spaces,
we prove that the problem is locally well-posed in $H^{-3/4}(\R)$ with $\beta <0$ and $\gamma>0.$
The new ingredient that we introduce in this paper is Lemmas 2.1-2.6.

\bigskip


\bigskip

\bigskip
\noindent {\bf Short Title:} The Cauchy problem for the  Ostrovsky Equation

\newpage

{\large\bf 1. Introduction}
\bigskip

\setcounter{Theorem}{0} \setcounter{Lemma}{0}

\setcounter{section}{1}

In this paper, we consider the Ostrovsky equation
$$
 \partial_{x}\left(u_{t}-\beta \partial_{x}^{3}u
  +\frac{1}{2}\partial_{x}(u^{2})\right)
    -\gamma u=0.
$$
This equation is a mathematical model of the propagation of weakly nonlinear long
waves in a rotating liquid. It was introduced by Ostrovsky in \cite{O} as a model
for weakly nonlinear long waves, by taking into account of the Coriolis force,
to describe the propagation of surface waves in the ocean in a rotating frame of reference.
The parameter $\gamma$ is a positive number and measures the effect of rotation,
and the parameter $\beta$ is a nozero real number of both signs and reflects
the type of dispersion of the media. When $\beta<0$, the equation has negative dispersion
and describes surface and internal waves in the ocean and surface waves
in a shallow channel with an uneven bottom. When $\beta>0$, the equation has
positive dispersion and describes capillary waves on the surface of
liquid or for oblique magneto-acoustic waves (see \cite{Be,GaSt,GiGrSt}).
In the absence of rotation (that is, $\gamma = 0$),
it becomes the Korteweg-de Vries equation. By changing variables the above
Ostrovsky equation can be written in the form
\begin{equation}
   u_{t}-\beta\partial_{x}^{3}u
    +\frac{1}{2}\partial_{x}(u^{2})- \gamma \partial_{x}^{-1} u
    =0.\label{1.01}
\end{equation}
The Ostrovsky equation has many important properties, such as solitary waves or soliton solutions,
etc., and it has closed relation to the KdV equation (see \cite{LL,LV,ZL,Tsu}).
It draws much attention of physists and mathematician.
Many people have investigated the Cauchy problem for (\ref{1.01}), for instance,
see \cite{GL, GH,HJ,HJ0, I,LL, IM1,IM2,IM3,IM4,IM5,LM,LV,VL,Z,ZL}.
By using the Fourier retriction norm method introduced in \cite{B,Bourgain-GAFA93},
 Isaza and Mej\'{\i}a \cite{IM1} proved that (\ref{1.01}) is locally well-posed
in $H^{s}(\R)$ with $s>-\frac{3}{4}$ in the negative dispersion case and
is locally well-posed in $H^{s}(\R)$ with $s>-\frac{1}{2}$ in the positive dispersion case.
Later they showed the ill-posedness in $H^{s}(\R)$ for  $s<-\frac{3}{4}$ (\cite{IM3}).
Recently, Tsugawa \cite{Tsu} proved the time local well-posedness in some anisotropic
Sobolev space $H^{s,a}(\R)$ with $s>-a/2-3/4$ and $0\leq a\leq 1$. The result
 includes the time local well-posedness in $H^{s}(\R)$  with $s>-3/4$ for
both positive and negative dispersion Ostrovsky equation. Thus, $s=-\frac{3}{4}$
is the critical regularity index for (\ref{1.01}) in the both dispersion cases.
Tsugawa considered also the weak rotation limit and proved that the solution
of the Ostrovsky equation converges to the solution of the KdV equation when
the rotation parameter $\gamma$ goes to 0. However, the well-posedness  of the Ostrovsky equation in the critical case
has been still open.

In this paper we study the Cauchy problem of the Ostrovsky equation (\ref{1.01}) with negative dispersion
complimented with the initial condition
\begin{equation}
u(0,x)=u_{0}(x),\quad x\in \R.\label{1.02}
\end{equation}
Compared with the KdV equation, the structure of the Ostrovsky equation is more complicated.
More precisely, the phase function of the KdV equation is the smooth function $\xi^{3}$, while
the phase function of the Ostrovsky equation is $\beta\xi^{3}+\frac{\gamma}{\xi}$,
which has a singular point $\xi=0$.  For the KdV equation, just as done in \cite{BT},
two simple identities $a^{3}+b^{3}-\frac{(a+b)^{3}}{4}=\frac{3}{4}(a+b)(a-b)^{2}$ and
$(a+b)^{3}-a^{3}-b^{3}=3ab(a+b)$ are valid to establish some key bilinear estimates,
which guarantee the wellposedness in the critical space $H^{s}(\R)$ with $s=-3/4$ (see
Guo \cite{G} and Kishimoto \cite{Kis}).
Obviously, if $u=u(x,t)$  is the solution to (\ref{1.01}), then $v(x,t)=\beta^{-1} u(x,-\beta^{-1} t)$
is the solution to the following equation
\begin{eqnarray}
v_{t}+ v_{xxx}+\frac{1}{2}(v^{2})_{x}+\beta^{-1} \gamma\partial_{x}^{-1}v=0\label{1.03}.
\end{eqnarray}
Without loss  of generality, throughout this paper, we can assume that $\beta=-1,\gamma=1$.
For the Ostrovsky equation in this paper,
the identities that we can utilize are
$$
 a^{3}-\frac{1}{a}+ b^{3}-\frac{1}{b}
    -\frac{(a+b)^{3}}{4}+\frac{4}{a+b}
  = \frac{3}{4}(a+b)(a-b)^{2}\left[1-\frac{4}{3 ab(a+b)^{2}}\right]
$$
and
$$
  (a+b)^{3}-\frac{1}{a+b}- a^{3}+\frac{1}{a}-b^{3}+\frac{1}{b}
  = 3ab(a+b)+ \frac{a^{2}+ab+b^{2}}{ab(a+b)}.
$$
These identities enable us to construct reasonable splitting of the spectral domains
so as to establish the crucial bilinear estimates for the local wellposedness of the problem.
As in \cite{BT,IK,IKT,T,Kis,G}, we may apply appropriate Besov-type spaces to establish the dyadic
bilinear estimates and finally we are able to show that
the Cauchy problem for (\ref{1.01}) is locally well-posed in $H^{-3/4}(\R)$ with $\beta <0,\gamma>0.$

We give some notations before stating the main results. Throughout this paper,
 $0<\epsilon<{10^{-4}}$.
 $C$ is a positive constant
 which may vary from line to line.  $A\sim B$ means that $|B|\leq |A|\leq 4|B|$.
 $A\gg B$ means that $|A|\geq 4|B|.$ $\psi(t)$ is a smooth function
 supported in $[-1,2]$ and equals to
 $1$ in $[0,1]$.  $\mathscr{F}$ denotes the Fourier transformation with respect to
 both space and time variables and $\mathscr{F}^{-1}$ denotes the inverse
 transformation of $\mathscr{F}$, while $\mathscr{F}_{x}$ denotes the Fourier transformation
 with respect to the space variable and $\mathscr{F}^{-1}_{x}$ denotes the
 inverse transformation
of $\mathscr{F}^{-1}_{x}$.
Denote
\begin{eqnarray*}
&&D:=\left\{(\tau,\xi)\in \R^{2}:|\xi|\leq1/8, |\tau|\geq |\xi|^{-3}\right\},\\
&&A_{j}=\left\{(\tau,\xi)\in \R^{2}: 2^{j}\leq \langle\xi\rangle<2^{j+1}\right\},\\
&&B_{k}=\left\{(\tau,\xi)\in \R^{2}: 2^{k}\leq \left\langle\tau-\xi^{3}+\frac{1}{\xi}
\right\rangle<2^{k+1}\right\},
\end{eqnarray*}
where $j,k$ are nonnegative integers. The restriction $|\xi|\leq \frac{1}{8}$
in the spectral domain $D$ is chosen according to the structure of the phase
function $-\xi^{3}+\frac{1}{\xi}$ of the Ostrovsky equation.

The Bourgain space
$
   X^{s, \>b}= \left\{u\in  \mathscr{S}^{'}(\R^{2})\, :\, \|u\|_{X^{s,\,b}}<\infty\right\}
$ is defined by the norm
$$
\|u\|_{X^{s, \>b}} =  \left\|\langle\xi\rangle^{s} \left\langle\tau-\xi^{3}
   + \frac{1}{\xi}\right\rangle^{b}\mathscr{F}u\right\|_{L_{\tau\xi}^{2}(\SR^{2})}.
$$
 The space
$
  X^{s,\,b,\> 1} = \left\{u\in  \mathscr{S}^{'}(\R^{2})\, :\, \|u\|^{X_{s,\>b,\>1}}<\infty\right\}
$  is defined by
where
\begin{eqnarray*}
    \|u\|_{X^{s,b,1}}
&=& \biggl\|\biggl(\left\|\langle\xi\rangle^{s}
    \left\langle\tau-\xi^{3}+\frac{1}{\xi}\right\rangle ^{b} \mathscr{F}u
     \right\|_{L_{\tau\xi}^{2}(A_{j}\cap B_{k})}\biggr)_{j,\>k\geq 0}
      \biggr\|_{\ell_{j}^{2}(\ell_{k}^{1})}
       \nonumber\\
&\sim& \biggl[\sum_{j}2^{js}\biggl(\sum_{k}2^{bk}\|\mathscr{F}u\|_{L_{\tau\xi}^{2}
        (A_{j}\cap B_{k})}\biggr)^{2}\biggr]^{1/2}.
\end{eqnarray*}
We shall also
use the norms $\|u\|_{X}$ and $\|u\|_{Y}$ of the spaces
$$
 X = \left\{u\in  \mathscr{S}^{'}(\R^{2})\, :\, \|u\|_{X}<\infty\right\},
 \,$$$$
 Y = \left\{u\in  \mathscr{S}^{'}(\R^{2})\, :\, \|u\|_{Y}<\infty\right\},
$$
where
\begin{eqnarray*}
    \|u\|_{X}
&=& \left\|\mathscr{F}^{-1}[\chi_{D^{c}}\mathscr{F}u]\right\|_{X^{-\frac{3}{4}, \frac{1}{2}, 1}}
    +\|\mathscr{F}^{-1}[\chi_{ D}\mathscr{F}u]\|_{X^{-\frac{3}{4}, \frac{1}{2}}},\\
    \|u\|_{Y}
&=& \left\|\langle\xi\rangle^{-3/4}\mathscr{F}u\right\|_{L_{\xi}^{2}L_{\tau}^{1}},
\end{eqnarray*}
where $D^{c}=\R^2_{\tau\xi}\setminus D$ is the complementary set of $D$ given above.
The spaces $\hat{X}, \hat{X}^{s,\>b,\>1}$ and $\hat{X}^{s,\>b}$
 are defined corresponding
to the following norms
\begin{eqnarray*}
\|f\|_{\hat{X}}=\|\mathscr{F}^{-1}f\|_{X},\quad
\|f\|_{\hat{X}^{s,b,1}}= \|\mathscr{F}^{-1}f\|_{X^{s,b,1}},
\quad
\|f\|_{\hat{X}^{s,b}}=\|\mathscr{F}^{-1}f\|_{X^{s,b}}.
\end{eqnarray*}
The space $ X_{T}$ is the restriction of $X$ onto the finite time interval $[-T,T]$ and
is defined according to the norm
 \begin{equation}
    \|u\|_{X_{T}} =\inf \left\{\|v\|_{X}:v\in X, u(t)=v(t)
 \>\> {\rm for} \> -T\leq t\leq T\right\}.\label{1.04}
 \end{equation}

The main result of this paper is as follow.

\begin{Theorem}\label{Thm1}
The Cauchy problem (\ref{1.01})(\ref{1.02}) is locally well-posed in $H^{-3/4}(\R)$ with
$\beta <0,\gamma >0$. That is,
for $u_{0} \in H^{-3/4}(\R)$, there exist a $T>0$ and a solution
$u\in C([-T, T]; H^{-3/4}(\R))$, and the solution map $u_0\mapsto u(t)$ is locally
Lipschitz continuous from $ H^{-3/4}(\R)$ into $C([-T, T]; H^{-3/4}(\R))$.
\end{Theorem}

The rest of the paper is arranged as follows. In Section 2,  we give some
preliminaries. In Section 3, we show two crucial dyadic bilinear estimates
and then apply them to establish bilinear estimates. In Section 4, we prove
the Theorem 1.1. Finally in Section 5, we give an appendix and show two examples
of the bilinear estimates in standard Bourgain spaces.

\noindent {\bf Remark:} Local well-posedness of  the  Cauchy problem for the Ostrovsky equation with positive dispersion  at the critical  regularity  and  the global well-posedness of  the Cauchy problem for the Ostrovsky
equation in $H^{s}(\R)$ with $s\geq -3/4$ has not been  established up to now, we will be devoted to
the problem later.

\bigskip
\bigskip
\bigskip

 \noindent{\large\bf 2. Preliminaries }

\setcounter{equation}{0}

\setcounter{Theorem}{0}

\setcounter{Lemma}{0}

\setcounter{section}{2}

In this section, we make some preparations. These includes the estimates for some convolutions
and basic inequality about the the phase functions which are used to get the dyadic bilinear
estimates in the Section 3. We also give some elementary estimates for the unitary group
corresponding to the Ostrovky equation.

\begin{Lemma}\label{Lemma2.1} Assume that $f,g\in \mathscr{S}'(\R^{2})$,
$\supp f\subset A_{j_{1}},\, \supp g \subset A_{j_{2}}$ and
\begin{eqnarray*}
K:=\inf\left\{|\xi_{1}-\xi_{2}|:\exists\> \tau_{1},\tau_{2},\> s.t.\>(\xi_{1},\tau_{1})\in
\supp f  ,\, (\xi_{2},\tau_{2})\in \supp g \right\}>0,
\end{eqnarray*}
If
\begin{eqnarray}
(\xi_{1},\tau_{1})\in
\supp f  ,\, (\xi_{2},\tau_{2})\in \supp g,\>\xi_{1}\xi_{2}<0\label{2.01}
\end{eqnarray}
or
\begin{eqnarray}
(\xi_{1},\tau_{1})\in
\supp f  ,\, (\xi_{2},\tau_{2})\in \supp g,\>\xi_{1}\xi_{2}\geq 0,\quad
\left|1-\frac{4}{3\xi^{2}\xi_{1}\xi_{2}}\right|
            > \frac{1}{2},\label{2.02}
\end{eqnarray}
then
\begin{eqnarray}
&&\||\xi|^{1/4}f\ast g\|_{L^{2}(\SR^{2})}\leq
C\|f\|_{\hat{X}^{0,\frac{1}{2},1}}
\|g\|_{\hat{X}^{0,\frac{1}{2},1}},\label{2.03}\\
&&\||\xi|^{1/2}f\ast g\|_{L^{2}(\SR^{2})}\leq CK^{-1/2}
\|f\|_{\hat{X}^{0,\frac{1}{2},1}}\|g\|_{\hat{X}^{0,\frac{1}{2},1}}.\label{2.04}
\end{eqnarray}
\end{Lemma}
{\bf Proof.} First we prove
\begin{eqnarray}
&&\left|\int_{\SR^{2}}\int_{\!\!\!\mbox{\scriptsize $
\begin{array}{l}
\xi=\xi_{1}+\xi_{2}\\
\tau=\tau_{1}+\tau_{2}
\end{array}
$}}(\chi_{B_{k_{1}}}f)(\xi_{1},\tau_{1})\>
    (\chi_{B_{k_{2}}}g)(\xi_{2},\tau_{2})\>
    |\xi|^{1/4} h(\tau,\xi)
     d\xi_{1}d\tau_{1}d\xi d\tau \right|\nonumber\\
&&\leq C2^{\frac{k_{1}+k_{2}}{2}}\|f\|_{L_{\xi\tau}^{2}(B_{k_{1}})}\|g\|_{L_{\xi\tau}^{2}(B_{k_{2}})}
\|h\|_{L_{\xi\tau}^{2}}.\label{2.05}
\end{eqnarray}
and
\begin{eqnarray}
&& \quad
    \left|\int_{\SR^{2}}\int_{\!\!\!\mbox{\scriptsize $
\begin{array}{l}
\xi=\xi_{1}+\xi_{2}\\
\tau=\tau_{1}+\tau_{2}
\end{array}
$}}
  (\chi_{B_{k_{1}}}f)(\xi_{1},\tau_{1})\>
  (\chi_{B_{k_{2}}}g)(\xi_{2},\tau_{2})\>
  |\xi|^{1/2}\, h(\xi,\tau)
  d\xi_{1}d\tau_{1}d\xi d\tau \right|\nonumber\\
&&
   \leq CK^{-1/2}2^{\frac{k_{1}+k_{2}}{2}}\>
  \|f\|_{L_{\xi\tau}^{2}(B_{k_{1}})}\>
   \|g\|_{L_{\xi\tau}^{2}(B_{k_{2}})}\>
   \|h\|_{L_{\xi\tau}^{2}}
   \label{2.06}
\end{eqnarray}
if (\ref{2.01}) or (\ref{2.02}) is valid.
By using the Cauchy-Schwartz inequality and the Fubini theorem, we obtain
\begin{eqnarray}
&& \quad
    \left|\int_{\SR^{2}}\int_{\!\!\!\mbox{\scriptsize $
\begin{array}{l}
\xi=\xi_{1}+\xi_{2}\\
\tau=\tau_{1}+\tau_{2}
\end{array}
$}}
    (\chi_{B_{k_{1}}}f)(\xi_{1},\tau_{1})\>
    (\chi_{B_{k_{2}}}g)(\xi_{2},\tau_{2})\>
    h(\xi,\tau)
    d\xi_{1}d\tau_{1}d\xi d\tau \right|\nonumber\\
&&
    \leq C\sup\limits_{(\xi,\>\tau)\in \SR^{2}}m_{1}(\xi,\tau)^{1/2}\>
    \|f\|_{L_{\xi\tau}^{2}(B_{k_{1}})}\>
    \|g\|_{L_{\xi\tau}^{2}(B_{k_{2}})}\>
    \|h\|_{L_{\xi\tau}^{2}},
    \label{2.07}
\end{eqnarray}
where
\begin{eqnarray*}
&\displaystyle
m_{1}(\tau,\xi)=\int\chi_{\Lambda_{1}}(\xi_{1},\tau_{1},\xi,\tau)d\xi_{1}d\tau_{1},
&\\
&\displaystyle
\Lambda_{1}:=\left\{(\xi_{1},\tau_{1},\xi,\tau)\in \R^{4}\,:\,
   (\xi_{1},\tau_{1})\in  \supp f,\>
   (\xi_{2},\tau_{2})\in  \supp g\right\},
&
\end{eqnarray*}
in which $\tau=\tau_1+\tau_2,\xi=\xi_1+\xi_2$ and (\ref{2.01}) or (\ref{2.02}) is valid. Thus, the proofs of (\ref{2.05}) and
(\ref{2.06}) are reduced to
\begin{eqnarray}
m_{1}(\tau,\xi)\leq C{\rm min}\left\{|\xi|^{-1/2}2^{k_{1}+k_{2}},K^{-1}|\xi|^{-1}
2^{k_{1}+k_{2}}\right\}.\label{2.08}
\end{eqnarray}
For fixed $\tau,\xi\neq 0$, let $E_1$ and $E_2$ be the projections of $\Lambda_{1}$ onto
the $\xi_1$-axis and $\tau_1$-axis respectively. We show
\begin{eqnarray}
&&\mes E_1 \le C \min\left\{ |\xi|^{-1/2}(2^{k_{1}/2}+2^{k_{2}/2} ),
              K^{-1} |\xi|^{-1} (2^{k_{1}}+2^{k_{2}})
      \right\},
      \label{2.09}\\
&&\mes E_2 \le C \min \left\{2^{k_{1}}, 2^{k_{2}}\right\},
      \label{2.010}
\end{eqnarray}
then (\ref{2.08}) follows.

As in the introduction, it is easily checked that
\begin{eqnarray}
 && \tau-\frac{\xi^{3}}{4}+\frac{4}{\xi}-\left(\tau_{1}-\xi_{1}^{3}
  + \frac{1}{\xi_{1}}\right)
  -\left(\tau_{2}-\xi_{2}^{3}+ \frac{1}{\xi_{2}}\right)\nonumber\\&&
= \frac{3}{4}\xi(\xi_{1}-\xi_{2})^{2}\left[1-\frac{4}{3\xi^{2}\xi_{1}\xi_{2}}\right].
\label{2.011}
\end{eqnarray}
From (\ref{2.011}),   we have that
\begin{equation}
\max \Biggl\{\frac{4|M-C(2^{k_{1}}+2^{k_{2}})|}
         {3|\xi|\left|1-\frac{4}{3\xi^{2}\xi_{1}\xi_{2}}\right|},
         K^{2}\Biggr\}
\leq  |2\xi_{1}-\xi|^{2}
\leq  \frac{4|M+C(2^{k_{1}}+2^{k_{2}})|}
         {3|\xi|\left|1-\frac{4}{3\xi^{2}\xi_{1}\xi_{2}}\right|},
\label{2.012}
\end{equation}
where
$
M=\tau-\frac{\xi^{3}}{4}+\frac{4}{\xi}
$ and $C$  is some generic positive constant.

\noindent Case (\ref{2.01}) holds:
in this case, $\xi_{1}\xi_{2}<0$.

\noindent  When
 $K\geq \Biggl\{
        \dfrac{4|M-C(2^{k_{1}}+2^{k_{2}})|}
              {3|\xi|\Bigl||1-\dfrac{4}{3\xi^{2}\xi_{1}\xi_{2}}\Bigr|}
       \Biggr\}^{1/2}$,
the length of the interval that
$|2\xi_{1}-\xi|$ lies in is bounded by
\begin{eqnarray}
&& \Biggl\{
    \frac{4|M+C(2^{k_{1}}+2^{k_{2}})|}
         {3|\xi|\left|1-\frac{4}{3\xi^{2}\xi_{1}\xi_{2}}\right|}
     \Biggr\}^{1/2}-K\nonumber\\
&&  \qquad
 =  \displaystyle
      \dfrac{\dfrac{4|M+C(2^{k_{1}}+2^{k_{2}})|}
                   {3|\xi|\left|1-\frac{4}{3\xi^{2}\xi_{1}\xi_{2}}\right|}
                   -K^{2}}
            {\Biggl\{\dfrac{4|M+C(2^{k_{1}}+2^{k_{2}})|}
                   {3|\xi|\left|1-\frac{4}{3\xi^{2}\xi_{1}\xi_{2}}\right|} \Biggr\}^{1/2}
                   +K}\nonumber\\
&&\qquad
   \leq
   \frac{
    \dfrac{4|M+C(2^{k_{1}}+2^{k_{2}})|}
          {3|\xi|\left|1-\frac{4}{3\xi^{2}\xi_{1}\xi_{2}}\right|}
  - \dfrac{4|M-C(2^{k_{1}}+2^{k_{2}})|}
          {3|\xi|\left|1-\frac{4}{3\xi^{2}\xi_{1}\xi_{2}}\right|}
        }
        {\Biggl\{\dfrac{4|M+C(2^{k_{1}}+2^{k_{2}})|}
                       {3|\xi|\left|1-\frac{4}{3\xi^{2}\xi_{1}\xi_{2}}\right|}
            \Biggr\}^{1/2}+K}\nonumber\\
&&\qquad\leq
      \dfrac{C(2^{k_{1}}+2^{k_{2}})}
            {|\xi|K\left|1-\frac{4}{3\xi^{2}\xi_{1}\xi_{2}}\right|}
\leq
      \dfrac{C(2^{k_{1}}+2^{k_{2}})}{|\xi|K}.
      \label{2.013}
\end{eqnarray}
From the first {\it inequality} of the above,
such length of the interval of $|2\xi_{1}-\xi|$ is also bounded by
\begin{eqnarray}
\frac{C(2^{k_{1}/2}+2^{k_{2}/2})}{|\xi|^{1/2}\left|1-\frac{4}
{3\xi^{2}\xi_{1}\xi_{2}}\right|^{1/2}}\leq
\frac{C(2^{k_{1}/2}+2^{k_{2}/2})}{|\xi|^{1/2}}.\label{2.014}
\end{eqnarray}
By (\ref{2.013}) and (\ref{2.014}), we obtain that the measure of $E_1$
in this part is bounded by
\begin{equation}
 C \min\left\{ |\xi|^{-1/2}(2^{k_{1}/2}+2^{k_{2}/2} ),
              K^{-1} |\xi|^{-1} (2^{k_{1}}+2^{k_{2}})
      \right\}.
      \label{2.015}
\end{equation}
When $K\leq \Biggl\{ \dfrac{4|M-C(2^{k_{1}}+2^{k_{2}})|}
{3|\xi|\left|1-\frac{4}{3\xi^{2}\xi_{1}\xi_{2}}\right|} \Biggr\}^{1/2}$,
the length of the interval of $|2\xi_{1}-\xi|$ is bounded by
\begin{eqnarray*}
&&\Biggl\{
    \frac{4|M+C(2^{k_{1}}+2^{k_{2}})|}
         {3|\xi|\left|1-\frac{4}{3\xi^{2}\xi_{1}\xi_{2}}\right|}
\Biggr\}^{1/2}
 -
\Biggl\{
    \frac{4|M-C(2^{k_{1}}+2^{k_{2}})|}
         {3|\xi|\left|1-\frac{4}{3\xi^{2}\xi_{1}\xi_{2}}\right|}
\Biggr\}^{1/2},\nonumber\\&&\leq  C
   \frac{
    \dfrac{4|M+C(2^{k_{1}}+2^{k_{2}})|}
          {3|\xi|\left|1-\frac{4}{3\xi^{2}\xi_{1}\xi_{2}}\right|}
  - \dfrac{4|M-C(2^{k_{1}}+2^{k_{2}})|}
          {3|\xi|\left|1-\frac{4}{3\xi^{2}\xi_{1}\xi_{2}}\right|}
        }
        {\Biggl\{\dfrac{4|M+C(2^{k_{1}}+2^{k_{2}})|}
                       {3|\xi|\left|1-\frac{4}{3\xi^{2}\xi_{1}\xi_{2}}\right|}
            \Biggr\}^{1/2}+K},\nonumber\\
\end{eqnarray*}
similar to (\ref{2.013}) and (\ref{2.014}), the measure of $E_1$
in this part is bounded by (\ref{2.015}).

\noindent Case (\ref{2.02}) holds: in this case, $\xi_{1}\xi_{2}\geq0$ and  $\left|1-\frac{4}
  {3\xi^{2}\xi_{1}\xi_{2}}\right| > \frac{1}{2}$, it can be proved similarly that
the measure of $E_1$ in this part is also bounded by (\ref{2.015}).
Recall that
$(\xi_1,\tau_1)\in B_{k_1}$ and  $(\xi_2,\tau_2)\in B_{k_2}$,
 \begin{eqnarray}
 \left|\tau_{1}-\xi_{1}^{3}+\frac{1}{\xi_{1}}\right|
 \leq C2^{k_{1}},
 \quad\left|\tau-\tau_{1}-(\xi-\xi_{1})^{3}+\frac{1}{\xi-\xi_{1}}\right|
 \leq
  C2^{k_{2}},\label{2.016}
 \end{eqnarray}
thus we get (\ref{2.010}).  Consequently, we have  (\ref{2.06}).

By using the Cauchy-Schwartz inequality and the triangle inequality, we have that
\begin{eqnarray}
&&\hspace{-1cm}\left|\int_{\SR^{2}}\int_{\!\!\!\mbox{\scriptsize $
\begin{array}{l}
\xi=\xi_{1}+\xi_{2}\\
\tau=\tau_{1}+\tau_{2}
\end{array}
$}}f(\xi_{1},\tau_{1})g(\xi_{2},\tau_{2})h(\xi,\tau)d\xi_{1}d\tau_{1}
d\xi d\tau \right|\nonumber\\
&&\hspace{-1cm}\leq C\sum_{k_{1}}\sum_{k_{2}}\left|\int_{\SR^{2}}
\int_{\!\!\!\mbox{\scriptsize $
\begin{array}{l}
\xi=\xi_{1}+\xi_{2}\\
\tau=\tau_{1}+\tau_{2}
\end{array}
$}}(\chi_{B_{k_{1}}}f)(\xi_{1},\tau_{1})(\chi_{B_{k_{2}}}g)(\xi_{2},\tau_{2})
h(\xi,\tau)d\xi_{1}d\tau_{1}d\xi d\tau \right|.\label{2.017}
\end{eqnarray}
Combining (\ref{2.07}), (\ref{2.08})  with (\ref{2.017}),
we have (\ref{2.03})-(\ref{2.04}).

We have completed the proof of Lemma 2.1.
\begin{Lemma}\label{Lemma2.2} Assume that $f,g\in \mathscr{S}'(\R)$,
$\supp f\subset A_{j_{1}},\, \supp g \subset A_{j_{2}}$.If
\begin{eqnarray}
 K_{1}:=\inf\left\{|\xi_{1}-\xi_{2}|:\exists\> \tau_{1},\tau_{2},\> s.t.\>(\xi_{1},\tau_{1})\in
\supp f  ,\, (\xi_{2},\tau_{2})\in \supp g\right\}\geq2,\label{2.018}
\end{eqnarray}
and
\begin{eqnarray}
 (\xi_{1},\tau_{1})\in
\supp f,\, (\xi_{2},\tau_{2})\in \supp g,\>\xi_{1}\xi_{2}\geq 0,\quad\left|1-\frac{4}{3\xi^{2}\xi_{1}\xi_{2}}\right|
            \leq \frac{1}{2},\label{2.019}
\end{eqnarray}
then
\begin{eqnarray}
&&\||\xi|^{1/4}f\ast g\|_{L^{2}(\SR^{2})}
\leq C\|f\|_{\hat{X}^{0,\frac{1}{2},1}}\|g\|_{\hat{X}^{0,\frac{1}{2},1}},\label{2.020}\\
&& \||\xi|^{1/2}f\ast g\|_{L^{2}(\SR^{2})}\leq CK_{1}^{-1/2}
\|f\|_{\hat{X}^{0,\frac{1}{2},1}}\|g\|_{\hat{X}^{0,\frac{1}{2},1}}.\label{2.021}
\end{eqnarray}
\end{Lemma}
{\bf Proof.}
First we prove
\begin{eqnarray}
&&\left|\int_{\SR^{2}}\int_{\!\!\!\mbox{\scriptsize $
\begin{array}{l}
\xi=\xi_{1}+\xi_{2}\\
\tau=\tau_{1}+\tau_{2}
\end{array}
$}}(\chi_{B_{k_{1}}}f)(\xi_{1},\tau_{1})\>
    (\chi_{B_{k_{2}}}g)(\xi_{2},\tau_{2})\>
    |\xi|^{1/4} h(\xi,\tau)
     d\xi_{1}d\tau_{1}d\xi d\tau \right|\nonumber\\
&&\leq C2^{\frac{k_{1}+k_{2}}{2}}\|f\|_{L_{\xi\tau}^{2}(B_{k_{1}})}\|g\|_{L_{\xi\tau}^{2}(B_{k_{2}})}
\|h\|_{L_{\xi\tau}^{2}}.\label{2.022}
\end{eqnarray}
and
\begin{eqnarray}
&& \quad
    \left|\int_{\SR^{2}}\int_{\!\!\!\mbox{\scriptsize $
\begin{array}{l}
\xi=\xi_{1}+\xi_{2}\\
\tau=\tau_{1}+\tau_{2}
\end{array}
$}}
  (\chi_{B_{k_{1}}}f)(\xi_{1},\tau_{1})\>
  (\chi_{B_{k_{2}}}g)(\xi_{2},\tau_{2})\>
  |\xi|^{1/2}\, h(\xi,\tau)
  d\xi_{1}d\tau_{1}d\xi d\tau \right|\nonumber\\
&&
   \leq CK_{1}^{-1/2}\>
   2^{\frac{k_{1}+k_{2}}{2}}\|f\|_{L_{\xi\tau}^{2}(B_{k_{1}})}\>
   \|g\|_{L_{\xi\tau}^{2}(B_{k_{2}})}\>
   \|h\|_{L_{\xi\tau}^{2}}
   \label{2.023}
\end{eqnarray}
if (\ref{2.019})  is valid.
By using the Cauchy-Schwartz inequality and the Fubini theorem, we obtain
\begin{eqnarray}
&& \quad
    \left|\int_{\SR^{2}}\int_{\!\!\!\mbox{\scriptsize $
\begin{array}{l}
\xi=\xi_{1}+\xi_{2}\\
\tau=\tau_{1}+\tau_{2}
\end{array}
$}}
    (\chi_{B_{k_{1}}}f)(\xi_{1},\tau_{1})\>
    (\chi_{B_{k_{2}}}g)(\xi_{2},\tau_{2})\>
    h(\xi,\tau)
    d\xi_{1}d\tau_{1}d\xi d\tau \right|\nonumber\\
&&
    \leq C\sup\limits_{(\xi,\tau\>)\in \SR^{2}}m_{2}(\xi,\tau)^{1/2}\>
    \|f\|_{L_{\xi\tau}^{2}(B_{k_{1}})}\>
    \|g\|_{L_{\xi\tau}^{2}(B_{k_{2}})}\>
    \|h\|_{L_{\xi\tau}^{2}},
    \label{2.024}
\end{eqnarray}
where
\begin{eqnarray*}
&\displaystyle
m_{2}(\xi,\tau)=\int\chi_{\Lambda_{2}}(\xi_{1},\tau_{1},\xi,\tau)d\xi_{1} d\tau_{1},
&\\
&\displaystyle
\Lambda_{2}:=\left\{(\xi_{1},\tau_{1},\xi,\tau)\in \R^{4}\,:\,
   (\xi_{1},\tau_{1})\in  \supp f,\>
   (\xi_{2},\tau_{2})\in \supp g\right\},
&
\end{eqnarray*}
in which $\tau=\tau_1+\tau_2,\xi=\xi_1+\xi_2$ and (\ref{2.019})  is valid. Thus, the proofs of (\ref{2.022}) and
(\ref{2.023}) are reduced to
\begin{eqnarray}
m_{2}(\tau,\xi)\leq C{\rm min}\left\{|\xi|^{-1/2}2^{k_{1}+k_{2}},K_{1}^{-1}|\xi|^{-1}
2^{k_{1}+k_{2}}\right\}.\label{2.025}
\end{eqnarray}
For fixed $\tau,\xi\neq 0$, let $E_3$ and $E_4$ be the projections of $\Lambda_{2}$ onto
the $\xi_1$-axis and $\tau_1$-axis respectively. We show
\begin{eqnarray}
&&\mes E_3 \le C \min\left\{ |\xi|^{-1/2}(2^{k_{1}/2}+2^{k_{2}/2} ),
              K_{1}^{-1} |\xi|^{-1} (2^{k_{1}}+2^{k_{2}})
      \right\},
      \label{2.026}\\
&&\mes E_4 \leq C \min \left\{2^{k_{1}}, 2^{k_{2}}\right\},
      \label{2.027}
\end{eqnarray}
then (\ref{2.025}) follows.
Since $\xi_{1}\xi_{2}\geq0$ and  $\left|1-\frac{4}{3
\xi^{2}\xi_{1}\xi_{2}}\right|\leq \frac{1}{2}$, we have
\begin{eqnarray*}
\frac{1}{2}\leq \frac{4}{3\xi^{2}\xi_{1}\xi_{2}}\leq\frac{3}{2}.
\end{eqnarray*}
The above inequality is equivalent to
\begin{eqnarray}
\frac{8}{9\xi^{2}}\leq \xi_{1}(\xi-\xi_{1})\leq \frac{8}{3\xi^{2}},\label{2.028}
\end{eqnarray}
from which we have
\begin{eqnarray}
     \frac12 \left({\xi-\sqrt{\xi^{2}-\frac{32}{9\xi^{2}}}}\,\right)
\leq \xi_{1}
\leq \frac12\left({\xi-\sqrt{\xi^{2}-\frac{32}{3\xi^{2}}}}\,\right)
     \label{2.029}
\end{eqnarray}
or
\begin{eqnarray}
     \frac12\left({\xi+\sqrt{\xi^{2}-\frac{32}{3\xi^{2}}}}\,\right)
\leq \xi_{1}
\leq \frac12\left({\xi+\sqrt{\xi^{2}-\frac{32}{9\xi^{2}}}}\,\right)
     \label{2.030}.
\end{eqnarray}
From (\ref{2.029})  and (\ref{2.030}),  we see that
the measure of $E_1$ in this part  is  bounded by
 \begin{eqnarray}
 C|\xi|^{-2}\leq C|\xi|^{-1}K_{1}^{-1}\leq C{\rm min}\left\{ \frac{(2^{k_{1}}+2^{k_{2}})}
 {|\xi|K_{1}},\frac{(2^{k_{1}/2}+2^{k_{2}/2})}{|\xi|^{1/2}}\right\}\label{2.031}
 \end{eqnarray}
 since $2\leq K_{1}\leq |\xi_{1}-\xi_{2}|\leq |\xi_{1}+\xi_{2}|=|\xi|.$
 Thus, we have (\ref{2.026}).
Since
$(\xi_1,\tau_1)\in B_{k_1}$ and  $(\xi_2,\tau_2)\in B_{k_2}$,
 \begin{eqnarray}
 \left|\tau_{1}-\xi_{1}^{3}+\frac{1}{\xi_{1}}\right|
 \leq C2^{k_{1}},
 \quad\left|\tau-\tau_{1}-(\xi-\xi_{1})^{3}+\frac{1}{\xi-\xi_{1}}\right|
 \leq
  C2^{k_{2}},\label{2.032}
 \end{eqnarray}
 we get (\ref{2.027}).  Consequently, we have  (\ref{2.025}).

By using the Cauchy-Schwartz inequality and the triangle inequality, we have that
\begin{eqnarray}
&&\hspace{-1cm}\left|\int_{\SR^{2}}\int_{\!\!\!\mbox{\scriptsize $
\begin{array}{l}
\xi=\xi_{1}+\xi_{2}\\
\tau=\tau_{1}+\tau_{2}
\end{array}
$}}f(\xi_{1},\tau_{1})g(\xi_{2},\tau_{2})h(\xi,\tau)d\xi_{1}d\tau_{1}
d\xi d\tau \right|\nonumber\\
&&\hspace{-1cm}\leq C\sum_{k_{1}}\sum_{k_{2}}\left|\int_{\SR^{2}}
\int_{\!\!\!\mbox{\scriptsize $
\begin{array}{l}
\xi=\xi_{1}+\xi_{2}\\
\tau=\tau_{1}+\tau_{2}
\end{array}
$}}(\chi_{B_{k_{1}}}f)(\xi_{1},\tau_{1})(\chi_{B_{k_{2}}}g)(\xi_{2},\tau_{2})
h(\xi,\tau)d\xi_{1}d\tau_{1}d\xi d\tau \right|.\label{2.033}
\end{eqnarray}
Combining (\ref{2.022}), (\ref{2.023})  with (\ref{2.033}),
we have (\ref{2.020})-(\ref{2.021}).

We have completed the proof of Lemma 2.2.

\noindent {\bf Remark 1:} From the proof process of (\ref{2.021}), to obtain (\ref{2.021}), it is sufficient to require that $K_{1}>0.$

\begin{Lemma}\label{Lemma2.3} Assume that $f,g\in \mathscr{S}'(\R^{2})$,
$\supp f\subset A_{j_{1}},\, \supp g \subset A_{j_{2}}$.
\begin{eqnarray*}
K_{2}:=\inf\left\{|\xi_{1}-\xi_{2}|:\exists\> \tau_{1},\tau_{2},\> s.t.\>(\xi_{1},\tau_{1})\in
\supp f  ,\, (\xi_{2},\tau_{2})\in \supp g\right\}>0,
\end{eqnarray*}
then
\begin{eqnarray*}
\||\xi|^{1/2}f\ast g\|_{L^{2}(\SR^{2})}\leq CK_{2}^{-1/2}
\|f\|_{\hat{X}^{0,\frac{1}{2},1}}\|g\|_{\hat{X}^{0,\frac{1}{2},1}}.
\end{eqnarray*}
\end{Lemma}
{\bf Proof.}  Combining Lemma 2.1 with Lemma 2.2 and Remark 1, we have that Lemma 2.3.
\begin{Lemma}\label{Lemma2.4} Assume that $f\in \mathscr{S}^{'}(\R^{2})$ , $g\in \mathscr{S}(\R^{2})$
with
$\mathop{\rm supp} f\subset A_{j}$ for some $j\geq0$ and $\,\Omega\subset \R^{2}$ has positive measure.
Let
\begin{eqnarray*}
K_{3}:=\inf\left\{|\xi_{1}+\xi|:\exists\, \tau,\tau_{1}\>
s.t.\> (\xi,\tau)\in \Omega,(\xi_{1},\tau_{1})\in \supp f
   \right\}>0.
\end{eqnarray*}
If
\begin{eqnarray}
(\xi,\tau)\in \Omega,(\xi_{1},\tau_{1})\in \supp f,\>\xi\xi_{1}>0\label{2.034}
\end{eqnarray}
or
\begin{eqnarray}
(\xi,\tau)\in \Omega  ,\, (\xi_{1},\tau_{1})\in \supp f,\>\xi\xi_{1}\geq 0,\quad
\left|1+\frac{4}{3\xi\xi_{1}\xi_{2}^{2}}\right|
            > \frac{1}{2}.\label{2.035}
\end{eqnarray}
Then, for any $k\geq0$, we have
\begin{eqnarray}
&&\|f\ast g\|_{L^{2}(B_{k})}
\leq
C2^{k/4}\|f\|_{\hat{X}^{0,\frac{1}{2},1}}\||\xi|^{-1/4}g\|_{L^{2}(\SR^{2})}
\label{2.036},\\
    &&\|f\ast g\|_{L^{2}(\Omega\cap B_{k})}
\leq  C2^{k/2}K_{3}^{-1/2} \|f\|_{\hat{X}_{0,\frac{1}{2},1}}\>
      \|\,|\xi|^{-1/2}g\|_{L^{2}(\SR^{2})}\>
      .
        \label{2.037}
\end{eqnarray}

\end{Lemma}
{\bf Proof.} First we prove
\begin{eqnarray}
&&\quad \left|\int_{\SR^{2}}\int_{\!\!\!\mbox{\scriptsize $
\begin{array}{l}
      \xi=\xi_{1}+\xi_{2}\\
      \tau=\tau_{1}+\tau_{2}
\end{array}
$}}
   (\chi_{B_{k_{1}}}f(\xi_{1},\tau_{1}))g(\xi_{2},\tau_{2})
   h(\xi,\tau)d\xi_{1}d\tau_{1}d\xi d\tau \right|\nonumber\\
&&
  \leq C2^{k_{1}/2}2^{k/4}\|f\|_{L_{\xi\tau}^{2}(B_{k_{1}})}\,\|\,|\xi|^{-1/4}g\|_{L_{\xi\tau}^{2}}
    \|h\|_{L_{\xi\tau}^{2}}\label{2.038}
\end{eqnarray}
for any $h\in L^{2}(\R^2)$ with $\mathop{\rm supp}h\subset B_{k}$  and
\begin{eqnarray}
&&\quad
  \left|\int_{\SR^{2}}\int_{\!\!\!\mbox{\scriptsize $
\begin{array}{l}
\xi=\xi_{1}+\xi_{2}\\
\tau=\tau_{1}+\tau_{2}
\end{array}
$}}(\chi_{B_{k_{1}}}f(\xi_{1},\tau_{1}))g(\xi_{2},\tau_{2})h(\xi,\tau)
d\xi_{1}d\tau_{1}d\xi d\tau \right|\nonumber\\
&&\leq CK_{3}^{-1/2}2^{\frac{k+k_{1}}{2}}\|f\|_{L_{\xi\tau}^{2}(B_{k_{1}})}\, \|\,|\xi|^{-1/2}g\|_{L_{\xi\tau}^{2}}
\|h\|_{L_{\xi\tau}^{2}}\label{2.039}
\end{eqnarray}
for any $h\in L^{2}(\R^2)$ with $\mathop{\rm supp}h\subset B_{k}\cap \Omega.$

By using the Cauchy-Schwartz inequality and the Fubini theorem, we obtain
\begin{eqnarray}
&&\quad
   \left|\int_{\SR^{2}}\int_{\!\!\!\mbox{\scriptsize $
\begin{array}{l}
\xi=\xi_{1}+\xi_{2}\\
\tau=\tau_{1}+\tau_{2}
\end{array}
$}}(\chi_{B_{k_{1}}}f(\xi_{1},\tau_{1}))g(\xi_{2},\tau_{2})h(\xi,\tau)
d\xi_{1}d\tau_{1}d\xi d\tau \right|
      \nonumber\\
&&    \leq C\sup\limits_{(\xi_{2}, \tau_{2})\in \SR^{2}}m_{3}(\xi_{2},\tau_{2})^{1/2}
     \|f\|_{L_{\xi\tau}^{2}}\|g\|_{L_{\xi\tau}^{2}}
     \|h\|_{L_{\xi\tau}^{2}},\label{2.040}
\end{eqnarray}
where
\begin{eqnarray*}
&\displaystyle
m_{3}(\xi_{2},\tau_{2})=\int\chi_{\Lambda_{2}}(\xi_{2},\tau_{2},\xi,\tau)d\xi d\tau,
&\\
&\displaystyle
\Lambda_{3}:=\left\{(\xi_{2},\tau_{2},\xi,\tau)\in \R^{4}\, :\,
(\xi_{1},\tau_{1})\in \supp f,\quad (\xi,\tau)\in \supp h\right\},
&
\end{eqnarray*}
in which $\tau=\tau_1+\tau_2,\xi=\xi_1+\xi_2$.  Hence, the proofs of
(\ref{2.038}) and (\ref{2.039}) are reduced to
\begin{eqnarray}
m_{3}(\xi_{2},\tau_{2})\leq C{\rm min}\left\{|\xi_{2}|^{-1/2}
2^{k/2+k_{1}},K_3^{-1}|\xi_{2}|^{-1}2^{k+k_{1}}\right\}.\label{2.041}
\end{eqnarray}
For fixed $\tau,\xi\neq 0$, $(\xi,\tau)\in B_{k}$, we let
$F_1$ and $F_2$ be the projections of $\Lambda_3$ onto
the the $\xi$-axis and $\tau$-axis respectively.
We shall show
\begin{eqnarray}
&&\mes F_1\le C\min \left\{|\xi_{2}|^{-1/2}(2^{k/2}+2^{k_{1}}),
             K_{3}^{-1}|\xi_{2}|^{-1}(2^{k}+2^{k_{1}})
      \right\},
       \label{2.042}\\
&&\mes F_2\le C\min \left\{2^{k},2^{k_{1}} \right\},
       \label{2.043}
\end{eqnarray}
then (\ref{2.041}) follows.

Similar to (\ref{2.011}), we have
\begin{eqnarray}
&&\tau_{2}-\frac{\xi_{2}^{3}}{4}+\frac{4}{\xi_{2}}
-\left(\tau-\xi^{3}+\frac{1}{\xi}\right)
+\left(\tau_{1}-\xi_{1}^{3}+\frac{1}{\xi_{1}}\right)
\nonumber\\&&=
\frac{3}{4}\xi_{2}(2\xi-\xi_{2})^{2}\left[1+\frac{4}
{3\xi\xi_{1}\xi_{2}^{2}}\right].\quad \label{2.044}
\end{eqnarray}
From (\ref{2.044}), we get
\begin{eqnarray}
{\max}\Biggl\{
           \frac{4|M_{1}-C(2^{k}+2^{k_{1}})|}
                {3|\xi_{2}|\left|1+\frac{4}{3\xi\xi_{1}\xi_{2}^{2}}
      \right|},\>  K_{3}^{2}\Biggr\}
\leq  |2\xi-\xi_{1}|^{2}
\leq  \frac{4|M_{1}+C(2^{k}+2^{k_{1}})|}
           {3|\xi_{2}|\left|1+\frac{4}{3\xi\xi_{1}\xi_{2}^{2}}\right|},
\label{2.045}
\end{eqnarray}
where
$
M_{1}=\tau_{2}-\frac{\xi_{2}^{3}}{4}+\frac{4}{\xi_{2}}.$

\noindent When (\ref{2.034}) holds:  in this case $\xi\xi_{1}>0$, thus,
we have that
\begin{eqnarray*}
\left|1+\frac{4}{3\xi\xi_{1}\xi_{2}^{2}}\right|>\frac{1}{2}.
\end{eqnarray*}
We consider
\begin{eqnarray*}
&&K_{3}\geq
    \Biggl\{
    \dfrac{4|M_{1}-C(2^{k}+2^{k_{1}})|}
          {3|\xi_{2}|\left|1+\frac{4}{3\xi\xi_{1}\xi_{2}^{2}}\right|}
    \Biggr\}^{1/2},\\
    &&K_{3}\leq \Biggl\{\dfrac{4|M_{1}-C(2^{k}+2^{k_{1}})|}
{3|\xi_{2}|\left|1+\frac{4}{3\xi\xi_{1}\xi_{2}^{2}}\right|}\Biggr\}^{1/2}.
\end{eqnarray*}
When
$K_{3}\geq
    \Biggl\{
    \dfrac{4|M_{1}-C(2^{k}+2^{k_{1}})|}
          {3|\xi_{2}|\left|1+\frac{4}{3\xi\xi_{1}\xi_{2}^{2}}\right|}
    \Biggr\}^{1/2}$,
from (\ref{2.045}), we have that the length of the interval that $|2\xi-\xi_{1}|$ lies in is bounded by
\begin{eqnarray}
&&
\Biggl\{\frac{4|M_{1}+C(2^{k}+2^{k_{1}})|}
             {3|\xi_{2}|\left|1+\frac{4}{3\xi\xi_{1}\xi_{2}^{2}}\right|}
 \Biggr\}^{1/2}-K_{3}\nonumber\\
&&\qquad
 =
 \dfrac{\dfrac{4|M_{1}+C (2^{k}+2^{k_{1}})|}
              {3|\xi_{2}|\left|1+\frac{4}{3\xi\xi_{1}\xi_{2}^{2}}\right|}
        -K_{3}^{2}}
      {\Biggl\{\dfrac{4|M_{1}+C(2^{k}+2^{k_{1}})|}
                     {3|\xi_{2}|\left|1+\frac{4}{3\xi\xi_{1}\xi_{2}^{2}}\right|}
       \Biggr\}^{1/2} + K_{3}}\nonumber\\
&&\qquad
   \leq \frac
{
    \dfrac{4|M_{1}+C(2^{k}+2^{k_{1}})|}
           {3|\xi_{2}|\left|1+\frac{4}{3\xi\xi_{1}\xi_{2}^{2}}\right|}
  - \dfrac{4|M_{1}-C(2^{k}+2^{k_{1}})|}
                {3|\xi_{2}|\left|1+\frac{4}{3\xi\xi_{1}\xi_{2}^{2}}\right|}
}
{
   \Biggl\{
    \dfrac{4|M_{1}+C(2^{k}+2^{k_{1}})|}
         {3|\xi_{2}|\left|1+\frac{4}{3\xi\xi_{1}\xi_{2}^{2}}\right|}
   \Biggr\}^{1/2} +K_{3}
}
\nonumber\\
&&\qquad
   \leq
    \dfrac{C(2^{k}+2^{k_{1}})}
         {|\xi_{2}|K_{3}\left|1+\frac{4}{3\xi\xi_{1}\xi_{2}^{2}}\right|}
   \leq \frac{C(2^{k}+2^{k_{1}})}{|\xi_{2}|K_{3}}.
    \label{2.046}
\end{eqnarray}
Moreover,  from the first {\it inequality} of the above, such length of the interval
of $|2\xi-\xi_1|$
 is also bounded by
\begin{equation}
   \frac{C(2^{k/2}+2^{k_{1}/2})}
        {|\xi_{2}|^{1/2}\left|1+\frac{4}{3\xi\xi_{1}\xi_{2}^{2}}
   \right|^{1/2}}
\leq
  \frac{C(2^{k/2}+2^{k_{1}/2})}{|\xi_{2}|^{1/2}}.
  \label{2.047}
\end{equation}
From (\ref{2.046}) and (\ref{2.047}), we infer that the measure of $F_1$ in this part
is bounded by
\begin{equation}
C\min \left\{|\xi_{2}|^{-1/2}(2^{k}+2^{k_{1}/2}),
             K_{3}^{-1}|\xi_{2}|^{-1}(2^{k}+2^{k_{1}})
      \right\}.
       \label{2.048}
\end{equation}
When $K_{3}\leq \Biggl\{\dfrac{4|M_{1}-C(2^{k}+2^{k_{1}})|}
{3|\xi_{2}|\left|1+\frac{4}{3\xi\xi_{1}\xi_{2}^{2}}\right|}\Biggr\}^{1/2}$,
the length of the interval of $|2\xi-\xi_{2}|$ is bounded by
\begin{eqnarray*}
&&\Biggl\{\frac{4|M_{1}+C(2^{k}+2^{k_{1}})|}
           {3|\xi_{2}|\left|1+\frac{4}{3\xi\xi_{1}\xi_{2}^{2}}\right|}
\Biggr\} ^{1/2}
  -
\Biggl\{
           \frac{4|M_{1}-C(2^{k}+2^{k_{1}})|}
                {3|\xi_{2}|\left|1+\frac{4}{3\xi\xi_{1}\xi_{2}^{2}}      \right|}
\Biggr\} ^{1/2}\nonumber\\&&\leq C\frac
{
    \dfrac{4|M_{1}+C(2^{k}+2^{k_{1}})|}
           {3|\xi_{2}|\left|1+\frac{4}{3\xi\xi_{1}\xi_{2}^{2}}\right|}
  - \dfrac{4|M_{1}-C(2^{k}+2^{k_{1}})|}
                {3|\xi_{2}|\left|1+\frac{4}{3\xi\xi_{1}\xi_{2}^{2}}\right|}
}
{
   \Biggl\{
    \dfrac{4|M_{1}+C(2^{k}+2^{k_{1}})|}
         {3|\xi_{2}|\left|1+\frac{4}{3\xi\xi_{1}\xi_{2}^{2}}\right|}
   \Biggr\}^{1/2} +K_{3}}.
\end{eqnarray*}
Similar to (\ref{2.046}) and (\ref{2.047}), the measure of $F_1$ in this part
is also bounded by (\ref{2.048}).

\noindent When (\ref{2.035}) holds: $\xi\xi_{1}\leq0$ and
$\left|1+\frac{4}{3\xi\xi_{1} \xi_{2}^{2}}\right| > \frac{1}{2}$,
it can be proved similarly that the measure of $F_1$ in this part
is also bounded by (\ref{2.048}).

 Recall that  $(\xi,\tau)\in B_{k}$, $(\xi_1,\tau_1)\in B_{k_1}$,
 \begin{eqnarray}
\left|\tau-\xi^{3}+ \frac{1}{\xi}\right|\leq C2^{k},
 \quad\left|\tau_{1}-\xi_{1}^{3}+ \frac{1}{\xi_{1}}\right|
\leq C2^{k_{1}},\label{2.049}
 \end{eqnarray}
we get the estimate (\ref{2.043}) for  $F_2$.
Consequently, we have  (\ref{2.041}).

By using the Cauchy-Schwartz inequality and the triangle inequality, we have that
\begin{eqnarray}
&&\hspace{-1cm}\left|\int_{\SR^{2}}\int_{\!\!\!\mbox{\scriptsize $
\begin{array}{l}
\xi=\xi_{1}+\xi_{2}\\
\tau=\tau_{1}+\tau_{2}
\end{array}
$}}f(\xi_{1},\tau_{1})g(\xi_{2},\tau_{2})h(\xi,\tau)
d\xi_{1}d\tau_{1}d\xi d\tau \right|\nonumber\\
&&\hspace{-1cm}\leq C\sum_{k_{1}}\left|\int_{\SR^{2}}
\int_{\!\!\!\mbox{\scriptsize $
\begin{array}{l}
\xi=\xi_{1}+\xi_{2}\\
\tau=\tau_{1}+\tau_{2}
\end{array}
$}}(\chi_{B_{k_{1}}}f(\xi_{1},\tau_{1}))g(\xi_{2},\tau_{2})
h(\xi,\tau)d\xi_{1}d\tau_{1}d\xi d\tau \right|.\label{2.050}
\end{eqnarray}
Combining (\ref{2.038}), (\ref{2.039})  with (\ref{2.050}),
we have (\ref{2.036}) and (\ref{2.037}).

We have completed the proof of Lemma \ref{Lemma2.4}.

\begin{Lemma}\label{Lemma2.5} Assume that $f\in \mathscr{S}^{'}(\R^{2})$ , $g\in \mathscr{S}(\R^{2})$
with
$\mathop{\rm supp} f\subset A_{j}$ for some $j\geq0$ and $\,\Omega\subset \R^{2}$ has positive measure.
Let
\begin{eqnarray*}
K_{4}:=\inf\left\{|\xi_{1}+\xi|:\exists\, \tau,\tau_{1}\>
s.t.\> (\xi,\tau)\in \Omega,(\xi_{1},\tau_{1})\in \supp f
   \right\}\geq 2.
\end{eqnarray*}
If
\begin{eqnarray}
(\xi,\tau)\in \Omega  ,\, (\xi_{1},\tau_{1})\in \supp f,\>\xi\xi_{1}\leq 0,\quad
\left|1+\frac{4}{3\xi\xi_{1}\xi_{2}^{2}}\right|
            \leq\frac{1}{2},\label{2.051}
\end{eqnarray}
then, for any $k\geq0$, we have
\begin{eqnarray}
&&\|f\ast g\|_{L^{2}(B_{k})}
\leq
C2^{k/4}\|f\|_{\hat{X}^{0,\frac{1}{2},1}}\||\xi|^{-1/4}g\|_{L^{2}(\SR^{2})}
\label{2.052},\\
    &&\|f\ast g\|_{L^{2}(\Omega\cap B_{k})}
\leq  C2^{k/2}K_{4}^{-1/2} \|f\|_{\hat{X}_{0,\frac{1}{2},1}}\>
      \|\,|\xi|^{-1/2}g\|_{L^{2}(\SR^{2})}\>
      .
        \label{2.053}
\end{eqnarray}
\end{Lemma}
{\bf Proof.} First we prove
\begin{eqnarray}
&&\quad \left|\int_{\SR^{2}}\int_{\!\!\!\mbox{\scriptsize $
\begin{array}{l}
      \xi=\xi_{1}+\xi_{2}\\
      \tau=\tau_{1}+\tau_{2}
\end{array}
$}}
   (\chi_{B_{k_{1}}}f(\xi_{1},\tau_{1}))g(\xi_{2},\tau_{2})
   h(\xi,\tau)d\xi_{1}d\tau_{1}d\xi d\tau \right|\nonumber\\
&&
  \leq C2^{k_{1}/2}2^{k/4}
   \|f\|_{L_{\xi\tau}^{2}(B_{k_{2}})}\,\|\,|\xi|^{-1/4}g\|_{L_{\xi\tau}^{2}}
    \|h\|_{L_{\xi\tau}^{2}}\label{2.054}
\end{eqnarray}
for any $h\in L^{2}(\R^2)$ with $\mathop{\rm supp}h\subset B_{k}$  and
\begin{eqnarray}
&&\quad
  \left|\int_{\SR^{2}}\int_{\!\!\!\mbox{\scriptsize $
\begin{array}{l}
\xi=\xi_{1}+\xi_{2}\\
\tau=\tau_{1}+\tau_{2}
\end{array}
$}}(\chi_{B_{k_{1}}}f(\xi_{1},\tau_{1}))g(\xi_{2},\tau_{2})h(\xi,\tau)
d\xi_{1}d\tau_{1}d\xi d\tau \right|\nonumber\\
&&\leq CK_{4}^{-1/2}2^{k/2+k_{1}/2}\|f\|_{L_{\xi\tau}^{2}(B_{k_{2}})}\, \|\,|\xi|^{-1/2}g\|_{L_{\xi\tau}^{2}}
\|h\|_{L_{\xi\tau}^{2}}\label{2.055}
\end{eqnarray}
for any $h\in L^{2}(\R^2)$ with $\mathop{\rm supp}h\subset B_{k}\cap \Omega.$

By using the Cauchy-Schwartz inequality and the Fubini theorem, we obtain
\begin{eqnarray}
&&\quad
   \left|\int_{\SR^{2}}\int_{\!\!\!\mbox{\scriptsize $
\begin{array}{l}
\xi=\xi_{1}+\xi_{2}\\
\tau=\tau_{1}+\tau_{2}
\end{array}
$}}(\chi_{B_{k_{1}}})f(\xi_{1},\tau_{1})g(\xi_{2},\tau_{2})h(\xi,\tau)
d\xi_{1}d\tau_{1}d\xi d\tau \right|
      \nonumber\\
&&    \leq C\sup\limits_{( \xi_{2}, \tau_{2})\in \SR^{2}}m_{4}(\xi_{2},\tau_{2})^{1/2}
     \|f\|_{L_{\xi\tau}^{2}}\|g\|_{L_{\xi\tau}^{2}}
     \|h\|_{L_{\xi\tau}^{2}},\label{2.056}
\end{eqnarray}
where
\begin{eqnarray*}
&\displaystyle
m_{4}(\xi_{2},\tau_{2})=\int\chi_{\Lambda_{4}}(\xi_{2},\tau_{2},\xi,\tau)d\xi d\tau,
&\\
&\displaystyle
\Lambda_{4}:=\left\{(\xi_{2},\tau_{2},\xi,\tau)\in \R^{4}\, :\,
(\xi_{1},\tau_{1})\in\supp f,(\xi,\tau)\in \supp g\right\},
&
\end{eqnarray*}
in which $\tau=\tau_1+\tau_2,\xi=\xi_1+\xi_2$.  Hence, the proofs of
(\ref{2.054}) and (\ref{2.055}) are reduced to
\begin{eqnarray}
m_{4}(\tau,\xi)\leq C{\rm min}\left\{|\xi_{2}|^{-1/2}
2^{k/2+k_{1}/2},K_4^{-1}|\xi_{2}|^{-1}2^{k+k_{1}}\right\}.\label{2.057}
\end{eqnarray}
For fixed $\tau,\xi\neq 0$, $(\xi,\tau)\in B_{k}$, we let
$F_3$ and $F_4$ be the projections of $\Lambda_4$ onto
the the $\xi_1$-axis and $\tau_1$-axis respectively.
We shall show
\begin{eqnarray}
&&\mes F_3\le C\min \left\{|\xi_{2}|^{-1/2}(2^{k/2}+2^{k_{1}/2}),
             K_{4}^{-1}|\xi_{2}|^{-1}(2^{k}+2^{k_{1}})
      \right\},
       \label{2.058}\\
&&\mes F_4\le C\min \left\{2^{k},2^{k_{1}} \right\},
       \label{2.059}
\end{eqnarray}
then (\ref{2.057}) follows.
When $\xi\xi_{1}\leq 0$ and
$\left|1+\frac{4}{3\xi\xi_{1}\xi_{2}^{2}}\right|\geq \frac{1}{2}$,
 we have
\begin{eqnarray}
\frac{8}{9\xi_{2}^{2}}\leq \xi(\xi-\xi_{2})\leq \frac{8}{3\xi_{2}^{2}}.\label{2.060}
\end{eqnarray}
From (\ref{2.060}), we have that
\begin{eqnarray}
\frac12\left({\xi_{2}+\sqrt{\xi_{2}^{2}+\frac{32}{9\xi_{2}^{2}}}}\,\right)
\leq \xi\leq
\frac12\left({\xi_{2}+\sqrt{\xi_{2}^{2}+\frac{32}{3\xi_{2}^{2}}}}\,\right)
\label{2.061}
\end{eqnarray}
or
\begin{eqnarray}
\frac12\left({\xi_{2}-\sqrt{\xi_{2}^{2}+\frac{32}{3\xi_{2}^{2}}}}\right)
\leq \xi
\leq \frac12\left({\xi_{2}-\sqrt{\xi_{2}^{2}+\frac{32}{9\xi_{2}^{2}}}}\,\right)
\label{2.062}.
\end{eqnarray}
From (\ref{2.061})  and (\ref{2.062}), we see that the measure of
 $F_3$  is  bounded by
 \begin{eqnarray*}
 C|\xi_{2}|^{-2}\leq C{\rm min}\left\{|\xi_{2}|^{-1/2}
2^{k/2+k_{1}/2},K_4^{-1}|\xi_{2}|^{-1}2^{k+k_{1}}\right\}.
 \end{eqnarray*}
Recall that  $(\xi,\tau)\in B_{k}$, $(\xi_1,\tau_1)\in B_{k_1}$,
 \begin{eqnarray}
\left|\tau-\xi^{3}+ \frac{1}{\xi}\right|\leq C2^{k},
 \quad\left|\tau_{1}-\xi_{1}^{3}+ \frac{1}{\xi_{1}}\right|
\leq C2^{k_{1}},\label{2.063}
 \end{eqnarray}
we get the estimate (\ref{2.059}) for  $F_4$.
Consequently, we have  (\ref{2.057}).

By using the Cauchy-Schwartz inequality and the triangle inequality, we have that
\begin{eqnarray}
&&\hspace{-1cm}\left|\int_{\SR^{2}}\int_{\!\!\!\mbox{\scriptsize $
\begin{array}{l}
\xi=\xi_{1}+\xi_{2}\\
\tau=\tau_{1}+\tau_{2}
\end{array}
$}}f(\xi_{1},\tau_{1})g(\xi_{2},\tau_{2})h(\xi,\tau)
d\xi_{1}d\tau_{1}d\xi d\tau \right|\nonumber\\
&&\hspace{-1cm}\leq C\sum_{k_{1}}\left|\int_{\SR^{2}}
\int_{\!\!\!\mbox{\scriptsize $
\begin{array}{l}
\xi=\xi_{1}+\xi_{2}\\
\tau=\tau_{1}+\tau_{2}
\end{array}
$}}(\chi_{B_{k_{1}}}f(\xi_{1},\tau_{1}))g(\xi_{2},\tau_{2})
h(\xi,\tau)d\xi_{1}d\tau_{1}d\xi d\tau \right|.\label{2.064}
\end{eqnarray}
Combining (\ref{2.054}), (\ref{2.055})  with (\ref{2.064}),
we have (\ref{2.052}) and (\ref{2.053}).

We have completed the proof of Lemma 2.5.

\noindent {\bf Remark 2:} From the proof process of (\ref{2.053}), to obtain (\ref{2.053}), it is sufficient to require that $K_{4}>0.$
\begin{Lemma}\label{Lemma2.6} Assume that $f\in \mathscr{S}^{'}(\R^{2})$ , $g\in \mathscr{S}(\R^{2})$
with
$\mathop{\rm supp} f\subset A_{j}$ for some $j\geq0$ and $\,\Omega\subset \R^{2}$ has positive measure.
If
\begin{eqnarray*}
K_{5}:=\inf\left\{|\xi_{1}+\xi|:\exists\, \tau,\tau_{1}\>
s.t.\> (\xi,\tau)\in \Omega,(\xi_{1},\tau_{1})\in \supp f
   \right\}>0.
\end{eqnarray*}
Then
\begin{eqnarray*}
\|f\ast g\|_{L^{2}(\Omega\cap B_{k})}
\leq  C2^{k/2}K_{5}^{-1/2} \|f\|_{\hat{X}_{0,\frac{1}{2},1}}\>
      \|\,|\xi|^{-1/2}g\|_{L^{2}(\SR^{2})}\>
\end{eqnarray*}
\end{Lemma}
{\bf Proof.} Combining Lemma 2.4 with Lemma 2.5, Remark 2, we have that Lemma 2.6.

\begin{Lemma}\label{Lemma2.7}
Let $j,N\in \N$, $\gamma_{0}=\frac 12 j\geq 2^{N+2},
\gamma_{n+1}=2\log_{2} {\gamma_{n}}$, $6\le \gamma_{_N}<8$. Then
\begin{eqnarray}
\sum_{n=0}^{N-1}\frac{1}{\sqrt{\gamma_{n}}}\label{2.065}
\end{eqnarray}
is bounded uniformly in $j$ and $N.$
\end{Lemma}
{\bf Proof.} We claim that
\begin{eqnarray}
\gamma_{n}\geq2^{N+2-n},\quad0\leq n\leq N-1.\label{2.066}
\end{eqnarray}
 Let
 $
 a_{n}=\gamma_{N-n},
 $
then (\ref{2.066}) is equivalent to
\begin{eqnarray}
 a_{n}=\gamma_{_{N-n}}\geq 2^{n+2}.\label{2.067}
 \end{eqnarray}
We prove (\ref{2.067}) by induction.

When $n=1$, $a_{1}=\gamma_{_{N-1}}=2^{\frac 12 \gamma_{_N}}\geq 2^{3}=8.$
We assume that for $n=k$, $a_{k}=\gamma_{_{N-k}}\geq 2^{k+2}.$ Then for $n=k+1,$
we have that
$
a_{k+1}= \gamma_{_{N-k-1}}=2^{\frac12{\gamma_{_{N-k}}}}
      \geq 2^{2^{k+1}}\geq  2^{k+3}.
$
Thus we have (\ref{2.067}). Consequently, we have
\begin{eqnarray}
\sum_{n=0}^{N-1}\frac{1}{\sqrt{\gamma_{n}}}\leq \sum_{n=0}^{N-1}\frac{1}{2^{(N+2-n)/2}}
=\frac{\sqrt{2}+1}{2}(1-2^{-N/2})\leq \frac{\sqrt{2}+1}{2}.\label{2.068}
\end{eqnarray}

We have completed the proof of Lemma \ref{Lemma2.7}.

\noindent{\bf Remark:} The conclusion of Lemma 2.7 can be found in page 460 of \cite{Kis}, however, the proof is not given.

\begin{Lemma}\label{Lemma2.8}
Let  $\tau=\tau_1+\tau_{2}, \xi=\xi_1+\xi_{2}.$ Then
\begin{eqnarray*}
&&\quad \max\left\{3|\xi\xi_{1}\xi_{2}|,
\frac{\xi_{1}^{2}+\xi_{1}\xi_{2}+\xi_{2}^{2}}{|\xi\xi_{1}\xi_{2}|}\right\}\\
&&
\leq\left|\left(\tau-\xi^{3}+\frac{1}{\xi}\right)-
\left(\tau_{1}-\xi_{1}^{3}+\frac{1}{\xi_{1}}\right)
-\left(\tau_{2}-\xi_{2}^{3}+\frac{1}{\xi_{2}}\right)\right|\nonumber\\
&&\leq 2{\rm max}\left\{3|\xi\xi_{1}\xi_{2}|,
\frac{\xi_{1}^{2}+\xi_{1}\xi_{2}+\xi_{2}^{2}}
{|\xi\xi_{1}\xi_{2}|}\right\}.
\end{eqnarray*}
\end{Lemma}
{\bf Proof.} By a direct computation, since $3\xi\xi_{1}\xi_{2}\times
\frac{\xi_{1}^{2}+\xi_{1}\xi_{2}+\xi_{2}^{2}}{\xi\xi_{1}\xi_{2}} >0,$ we have that
\begin{eqnarray}
&&\quad \max\left\{3|\xi\xi_{1}\xi_{2}|,
            \frac{\xi_{1}^{2}+\xi_{1}\xi_{2}+\xi_{2}^{2}}
            {|\xi\xi_{1}\xi_{2}|}\right\}
            \nonumber\\
&&\leq\left|\left(\tau-\xi^{3}+\frac{1}{\xi}\right)-
       \left(\tau_{1}-\xi_{1}^{3}+\frac{1}{\xi_{1}}\right)
        -\left(\tau_{2}-\xi_{2}^{3}+\frac{1}{\xi_{2}}\right)\right|
          \nonumber\\
&&
 =\left|3\xi\xi_{1}\xi_{2}+\frac{\xi_{1}^{2}
   +\xi_{1}\xi_{2}+\xi_{2}^{2}}{\xi\xi_{1}\xi_{2}}\right|
    \nonumber\\
&&\leq 2\max\left\{3|\xi\xi_{1}\xi_{2}|,
         \frac{\xi_{1}^{2}+\xi_{1}\xi_{2}+\xi_{2}^{2}}
              {|\xi\xi_{1}\xi_{2}|}\right\}
          \label{2.069}.
\end{eqnarray}

We have completed the proof of Lemma \ref{Lemma2.8}.

\begin{Lemma}\label{Lemma2.9}  The space $\hat{X}$ has the following properties.
(i) For any $b>1/2$, there exists $C>0$ such that
\begin{eqnarray}
\||f\|_{\hat{X}}&\leq&
C\|f\|_{\hat{X}^{-3/4,\>b}}.\label{2.070}
\end{eqnarray}
(ii) For $1< p\leq 2$, there exists $C>0$ such that
\begin{eqnarray}
&&\|\langle\xi \rangle^{-3/4} f\|_{L_{\xi}^{2}L_{\tau}^{p}}\leq
C\|f\|_{\hat{X}},\label{2.071}\\
&&\|\langle\xi \rangle^{-3/4} f\|_{L_{\xi}^{2}L_{\tau}^{1}}\leq
C\|f\|_{\hat{X}^{-3/4,\frac{1}{2},1}}\label{2.072}.
\end{eqnarray}

\end{Lemma}
{\bf Proof.} (i) can be proved similarly to (i) of \cite{Kis}. (\ref{2.071})
can be proved similarly to $1<p\leq 2$ of  (ii) in \cite{Kis}.
(\ref{2.072})  can be proved similarly to \cite{BT}.

We have completed the proof of Lemma \ref{Lemma2.9}.

\begin{Lemma}\label{Lemma2.10}  Let $e^{-t(-\partial_{x}^{3}+\partial_{x}^{-1})}u_{0}$
be the solution to the linear equation (\ref{1.01}). Then we have the following estimates
\begin{eqnarray}
&&\hspace{-1cm}\left\|e^{-t(-\partial_{x}^{3}+\partial_{x}^{-1})}u_{0}\right\|_{X_{T}}+\sup
\limits_{-T\leq t\leq T}\|e^{-t(-\partial_{x}^{3}+\partial_{x}^{-1})}u_{0}\|_{H_{x}^{-3/4}(\SR)}\leq
C\|u_{0}\|_{H_{x}^{-3/4}(\SR)}\label{2.073}
\end{eqnarray}
and
\begin{eqnarray}
&&\left\|\int_{0}^{t}e^{-(t-s)(-\partial_{x}^{3}+\partial_{x}^{-1})}
         F(s)ds\right\|_{X_{T}}%\nonumber\\
%&&\qquad
       +\sup\limits_{-T\leq t\leq T}\left\|\int_{0}^{t}
         e^{-(t-s)(-\partial_{x}^{3}+\partial_{x}^{-1})}
           F(s)ds\right\|_{H_{x}^{-3/4}(\SR)}
            \qquad
             \nonumber\\
&&\hspace{-5mm}\leq
 C\left\|\mathscr{F}^{-1}
   \Bigl(\left\langle\tau-\xi^{3}
   +\frac{1}{\xi}\right\rangle^{-1}\mathscr{F}F\Bigr)\right\|_{X}
%     \nonumber\\
%&& \qquad {}
   +\left\|\mathscr{F}^{-1}\Bigl(\left\langle\tau-\xi^{3}
    +\frac{1}{\xi}\right\rangle^{-1}\mathscr{F}F\Bigr)\right\|_{Y},
      \label{2.074}
\end{eqnarray}
where $0\leq T\leq1.$
\end{Lemma}

Lemma \ref{Lemma2.10} can be proved similarly to Lemma 4.1 of \cite{Kis}.

\bigskip
\bigskip

\noindent{\large\bf 3. Bilinear estimates }

\setcounter{equation}{0}

 \setcounter{Theorem}{0}

\setcounter{Lemma}{0}

 \setcounter{section}{3}

In this section, we give the proof of Lemmas 3.1-3.2 which is the core of this paper.

\begin{Lemma}\label{Lemma3.1}  Suppose $f, g \in \mathscr{S}'(\R^{2})$,
$\supp f\subset A_{j_{1}}$ and $\supp g\subset A_{j_{2}}$. Then we have
\begin{eqnarray}
     \left\|I_{A_{j}}\left\langle\tau-\xi^{3}+\frac {1}{\xi}\right\rangle^{-1}
      \xi f*g\right\|_{\hat{X}}
&\leq&  C(j,j_{1},j_{2})\|f\|_{\hat{X}}\|g\|_{\hat{X}},
        \label{3.01}\\
     \left\|I_{A_{j}}\langle\xi\rangle^{-3/4}
      \left\langle\tau-\xi^{3}+\frac{1}{\xi}\right\rangle^{-1}
       \xi f*g\right\|_{L_{\xi}^{2}L_{\tau}^{1}}
&\leq& C(j,j_{1},j_{2})\|f\|_{\hat{X}}\|g\|_{\hat{X}}
         \label{3.02}
\end{eqnarray}
for $j\geq0$ in the following cases.

(i)  At least two of $j,j_{1},j_{2}$ are less than 30 and $C(j,j_{1},j_{2})\sim 1$.

(ii) $j_{1},j_{2}\geq30$,  $|j_{1}-j_{2}|\leq 10$,   $0<j<j_{1}-9$ and $C(j,j_{1},j_{2})\sim 2^{-\frac{3}{8} j}$.

(iii) $j,j_{1}\geq 30,$ $|j-j_{1}|\leq 10$, $0<j_{2}<j-10$  and  $C(j,j_{1},j_{2})\sim2^{-\frac{1}{4}(j-j_{2})}.$

(iv) $j,j_{2}\geq 30,$ $|j-j_{2}|\leq 10$, $0<j_{1}<j-10$ and $C(j,j_{1},j_{2})\sim2^{-\frac{1}{4}(j-j_{1})}.$

(v) $j,j_{1},j_{2}\geq 30$, $|j-j_{1}|\leq 10$, $|j-j_{2}|\leq 10$ and $C(j,j_{1},j_{2})\sim 1.$

(vi) $j_{1},j_{2}\geq30$,   $j=0$ and $C(j,j_{1},j_{2})\sim 1.$

(vii)  $j,j_{1}\geq 30$,   $j_{2}=0$ and $C(j,j_{1},j_{2})\sim 1.$

(viii)  $j,j_{2}\geq 30$,  $j_{1}=0,$ and $C(j,j_{1},j_{2})\sim 1.$

\end{Lemma}

{\bf Proof.} (i)  In this case we may assume that $j,j_{1},j_{2}$ are all less than 40.
By using   the Young inequality and (\ref{2.071})-(\ref{2.072}), we have
\begin{eqnarray}
&&\left\|I_{A_{j}}\left\langle\tau-\xi^{3}+\frac{1}{\xi}\right\rangle^{-1}\xi f*g\right\|_{\hat{X}}\leq C
\left\|I_{A_{j}}\left\langle\tau-\xi^{3}+\frac{1}{\xi}\right\rangle^{-1}\xi f*g\right\|_{\hat{X}^{-3/4,1/2,1}}
\nonumber\\&&\leq C\|f*g\|_{L^{2}}\leq C
\|f\|_{L_{\xi}^{2}L_{\tau}^{4/3}}\|g\|_{L_{\xi}^{2}L_{\tau}^{4/3}}\leq C
\|f\|_{\hat{X}}\|g\|_{\hat{X}}\label{3.03}
\end{eqnarray}
and
\begin{eqnarray}
&&\left\|I_{A_{j}}\left\langle\tau-\xi^{3}+\frac{1}{\xi}\right\rangle^{-1}\langle\xi\rangle^{-3/4}\xi f*g\right\|_{L_{\xi}^{2}L_{\tau}^{1}}
\leq C\left\| f*g\right\|_{L_{\xi}^{2}L_{\tau}^{2}}\nonumber\\&&\leq C
\|f\|_{L_{\xi}^{2}L_{\tau}^{4/3}}\|g\|_{L_{\xi}^{2}L_{\tau}^{4/3}}\leq C\|f\|_{\hat{X}}
\|g\|_{\hat{X}}\label{3.04}.
\end{eqnarray}
(ii) In this case,
we restrict  $f$ to $B_{k_{1}}$ and $g$ to $B_{k_{2}}$, by using Lemma 2.8,
we have that
\begin{eqnarray}
2^{k_{max}}:=2^{{\rm max}\{k,\>k_{1},\>k_{2}\}}\geq C2^{j+2j_{1}}\label{3.05}.
\end{eqnarray}
When $2^{k}\geq C2^{j+2j_{1}}$ which yields  that $2^{j/4}2^{-k/2}\leq C2^{-3j/4}2^{-j_{1}}2^{j/2},$
by using Lemma 2.3, we have that
\begin{eqnarray}
&&\left\|I_{A_{j}}\left\langle\tau-\xi^{3}+\frac{1}{\xi}\right\rangle^{-1}\xi
f*g\right\|_{\hat{X}}\leq C2^{j/4}
\sum_{k\geq0}2^{-k/2}\left\| f*g\right\|_{L^{2}(A_{j}\cap B_{k})}\nonumber\\
&&\leq C2^{-\frac{3}{4}j}2^{-j_{1}}\sum_{k=j+2j_{1}+O(1)}\left\||\xi|^{1/2}
f*g\right\|_{L^{2}}\leq C2^{-\frac{3}{4}j}2^{-\frac{3}{2}j_{1}}
\|f\|_{\hat{X}^{0,1/2,1}}\|g\|_{\hat{X}^{0,1/2,1}}\nonumber\\&&
\leq C2^{-\frac{3}{4}j}
\|f\|_{\hat{X}^{-\frac{3}{4},1/2,1}}\|g\|_{\hat{X}^{-\frac{3}{4},1/2,1}}\nonumber\\&&\leq C2^{-\frac{3}{8}j}\|f\|_{\hat{X}}
\|g\|_{\hat{X}}
\label{3.06}.
\end{eqnarray}
When $2^{k_{1}}\geq C2^{j+2j_{1}}$ which yields that $2^{-\frac{3}{8}(j+2j_{1})}2^{3k_{1}/8}
2^{\frac{1}{8}(k_{1}-k)}\geq C,$  by using Lemma 2.6, we have that
\begin{eqnarray}
&&\hspace{-0.5cm}\left\|I_{A_{j}}\left\langle\tau-\xi^{3}+
\frac{1}{\xi}\right\rangle^{-1}\xi f*g\right\|_{\hat{X}}\nonumber\\&&\leq C2^{-j/8}
2^{-3j_{1}/4}\sum_{k\geq0}2^{-5k/8}\left\|\left(\left\langle\tau-\xi^{3}
+\frac{1}{\xi}\right\rangle^{1/2} f\right)*g
\right\|_{L^{2}(A_{j}\cap B_{k})}\nonumber\\
&&\leq C2^{-\frac{1}{8}j}2^{-\frac{7}{4}j_{1}}\|f\|_{\hat{X}^{0,1/2,1}}
\|g\|_{\hat{X}^{0,1/2,1}}\nonumber\\&&
\leq C2^{-\frac{3}{8}j}
\|f\|_{\hat{X}^{-\frac{3}{4},1/2,1}}\|g\|_{\hat{X}^{-\frac{3}{4},1/2,1}}\nonumber\\&&\leq C2^{-\frac{3}{8}j}\|f\|_{\hat{X}}
\|g\|_{\hat{X}}\label{3.07}.
\end{eqnarray}
When $2^{k_{2}}\geq C2^{j+2j_{1}}$,  this case can be treated similarly
to case $2^{k_{1}}\geq C2^{j+2j_{1}}$.

In this case, by using the Cauchy-Schwartz inequality  with respect to $\tau$ and a proof similar to the above cases, we have that
\begin{eqnarray}
&&\left\|I_{A_{j}}\left\langle\tau-\xi^{3}+\frac{1}{\xi}\right\rangle^{-1}
\xi \langle \xi\rangle ^{-3/4}f*g\right\|_{L_{\xi}^{2}L_{\tau}^{1}}
\nonumber\\&&\leq C2^{j/4}\sum_{k\geq0}2^{-k/2}\left\| f*g\right\|
_{L_{\xi}^{2}L_{\tau}^{2}}\nonumber\\&&
\leq C\left[2^{-\frac{3}{4}j}2^{-\frac{3}{2}j_{1}}
+2^{-\frac{1}{8}j}2^{-\frac{7}{4}j_{1}}\right]\|f\|_{\hat{X}^{0,1/2,1}}
\|g\|_{\hat{X}^{0,1/2,1}}\nonumber\\&&
\leq C2^{-\frac{3}{8}j}
\|f\|_{\hat{X}^{-\frac{3}{4},1/2,1}}\|g\|_{\hat{X}^{-\frac{3}{4},1/2,1}}\nonumber\\&&\leq C2^{-\frac{3}{8}j}\|f\|_{\hat{X}}
\|g\|_{\hat{X}}\label{3.08}.
\end{eqnarray}
(iii)
In this case, from Lemma \ref{Lemma2.8}, we have that $2^{k_{\rm max}}\geq C|\xi
\xi_{1}\xi_{2}|\geq C2^{2j+j_{2}}$

In this case, the left hand side of (\ref{3.01})-(\ref{3.02}) can be bounded by
\begin{eqnarray}
C2^{j/4}\sum_{k\geq0}2^{-k/2}\left\| f*g\right\|_{L_{\xi}^{2}L_{\tau}^{2}}\label{3.09}.
\end{eqnarray}
When $2^{k}\sim 2^{k_{\rm max}}\geq C2^{2j+j_{2}}$, with the aid of Lemma 2.3,  (\ref{3.09})
can be bounded by
\begin{eqnarray}
&&C2^{-\frac{7}{4}j}2^{-j_{2}/2}\|f\|_{\hat{X}^{0,1/2,1}}\|g\|_{\hat{X}^{0,1/2,1}}\nonumber\\&&\leq C2^{-j+\frac{j_{2}}{4}}\|f\|_{\hat{X}^{-\frac{3}{4},1/2,1}}\|g\|_{\hat{X}^{-\frac{3}{4},1/2,1}}\nonumber\\
&&\leq C2^{-j+\frac{j_{2}}{4}}\|f\|_{\hat{X}}
\|g\|_{\hat{X}}\label{3.010}.
\end{eqnarray}
When $2^{k_{1}}\sim 2^{k_{\rm max}}\geq C2^{2j+j_{2}}$, with the aid of (\ref{2.036}) and the fact
that $2^{k_{1}/2}2^{-\frac{1}{2}(2j+j_{2})}\geq C,$ we have that (\ref{3.09})  can be bounded by
\begin{eqnarray}
&&C2^{-3j/4}2^{-j_{2}/2}2^{k_{1}/2}\sum_{k\geq0}2^{-k/2}\left\| f*g\right\|_{L_{\xi}^{2}L_{\tau}^{2}}\nonumber\\&&
\leq C2^{-3j/4}2^{-j_{2}/2}\sum_{k\geq 0}2^{-k/4}2^{k_{1}/2}\||\xi|^{-1/4}f\|_{L^{2}}\|g\|_{\hat{X}^{0,1/2,1}}\nonumber\\&&\leq C2^{-\frac{1}{4}(j-j_{2})}\|f\|_{\hat{X}^{-\frac{3}{4},1/2,1}}\|g\|_{\hat{X}^{-\frac{3}{4},1/2,1}}\nonumber\\
&&\leq C2^{-\frac{1}{4}(j-j_{2})}\|f\|_{\hat{X}}
\|g\|_{\hat{X}}\label{3.011}.
\end{eqnarray}
When $2^{k_{2}}\sim 2^{k_{\rm max}}\geq C2^{2j+j_{2}}$, by using the fact
that $2^{k_{2}/2}2^{-\frac{1}{2}(2j+j_{2})}\geq C,$ in this case $2^{j_{2}}\ll 2^{j}$, we have that
(\ref{3.09})  can be bounded by
\begin{eqnarray}
C2^{-3j/4}2^{-j_{2}/2}2^{k_{2}/2}\sum_{k\geq0}2^{-k/2}\left\| f*g\right\|_{L_{\xi}^{2}L_{\tau}^{2}}\label{3.012},
\end{eqnarray}
if $k\leq 10j,$ by using  (\ref{2.037}) and $2^{j/4}\geq j(j\geq 20)$, we have that (\ref{3.012})  can be bounded by
\begin{eqnarray}
&&C2^{-5j/4}2^{-j_{2}}\sum_{k\geq0}2^{-k/4}\|f\|_{\hat{X}^{0,1/2,1}}\|g\|_{\hat{X}^{0,1/2,1}}\nonumber\\&&\leq Cj2^{-\frac{1}{2}j-\frac{j_{2}}{4}}\|f\|_{\hat{X}^{-\frac{3}{4},1/2,1}}\|g\|_{\hat{X}^{-\frac{3}{4},1/2,1}}\nonumber\\&&\leq C2^{-\frac{1}{4}(j-j_{2})}\|f\|_{\hat{X}^{-\frac{3}{4},1/2,1}}\|g\|_{\hat{X}^{-\frac{3}{4},1/2,1}}\nonumber\\
&&\leq C2^{-\frac{1}{4}(j-j_{2})}\|f\|_{\hat{X}}
\|g\|_{\hat{X}}\label{3.013},
\end{eqnarray}
if $k\geq 10j,$
we have that (\ref{3.012})  can be bounded by
\begin{eqnarray}
&&C2^{-5j/4}2^{-j_{2}}\sum_{k\geq0}2^{-k/4}\|f\|_{\hat{X}^{0,1/2,1}}
\|g\|_{\hat{X}^{0,1/2,1}}\nonumber\\&&\leq
C2^{-15j/4}2^{-j_{2}}\|f\|_{\hat{X}^{0,1/2,1}}\|g\|_{\hat{X}^{0,1/2,1}}
\nonumber\\&&\leq C2^{-3j-j_{2}/4}
\|f\|_{\hat{X}^{-\frac{3}{4},1/2,1}}\|g\|_{\hat{X}^{-\frac{3}{4},1/2,1}}
\nonumber\\
&&\leq C2^{-3j-j_{2}/4}\|f\|_{\hat{X}}
\|g\|_{\hat{X}}\label{3.014}.
\end{eqnarray}
(iv) This case can be proved similarly to case (iii).

\noindent(v) In this case, from Lemma \ref{Lemma2.8}, we have that $2^{k_{\rm max}}\geq C|\xi
\xi_{1}\xi_{2}|\geq C2^{3j}\sim 2^{3j_{1}}\sim 2^{3j_{2}}.$
In this case, the left hand side of (\ref{3.01})-(\ref{3.02}) can be bounded by
\begin{eqnarray}
C2^{j/4}\sum_{k\geq0}2^{-k/2}\left\| f*g\right\|_{L_{\xi}^{2}L_{\tau}^{2}}\label{3.015}.
\end{eqnarray}
When $2^{k}\sim 2^{k_{\rm max}}\geq C2^{2j+j_{2}}$, with the aid of (\ref{2.03}),  (\ref{3.015})
can be bounded by
\begin{eqnarray}
&&C2^{-3j/2}\|f\|_{\hat{X}^{0,1/2,1}}\|g\|_{\hat{X}^{0,1/2,1}}
\leq C\|f\|_{\hat{X}^{-\frac{3}{4},1/2,1}}\|g\|_{\hat{X}^{-\frac{3}{4},1/2,1}}\nonumber\\
&&\leq C\|f\|_{\hat{X}}
\|g\|_{\hat{X}}\label{3.016}.
\end{eqnarray}
When $2^{k_{1}}\sim 2^{k_{\rm max}}\geq C2^{3j}$, with the aid of (\ref{2.036}) and the fact
that $2^{k_{1}/2}2^{-\frac{3}{2}j}\geq C,$ (\ref{3.015})  can be bounded by
\begin{eqnarray}
&&C2^{-5j/4}2^{k_{1}/2}\sum_{k\geq0}2^{-k/2}\left\| f*g\right\|_{L_{\xi}^{2}
L_{\tau}^{2}}\nonumber\\&&
\leq C2^{-\frac{5}{4}j}\|f\|_{\hat{X}^{0,1/2,1}}\|g\|_{\hat{X}^{0,1/2,1}}\nonumber\\&&\leq C
\|f\|_{\hat{X}^{-\frac{3}{4},1/2,1}}\|g\|_{\hat{X}^{-\frac{3}{4},1/2,1}}\nonumber\\&&
\leq C\|f\|_{\hat{X}}
\|g\|_{\hat{X}}\label{3.017}.
\end{eqnarray}
When $2^{k_{2}}\sim 2^{k_{\rm max}}\geq C2^{2j+j_{2}}$, this case can be proved similarly to case
$2^{k_{1}}\sim 2^{k_{\rm max}}\geq C2^{3j}$.

\noindent (vi)  (a) When $2^{k_{1}}\sim 2^{k_{\rm max}}\geq 2^{2j_{1}+j}$, this case can be proved
similarly to $2^{k_{1}}\sim 2^{k_{\rm max}}$ of (ii).

(b) When $2^{k_{2}}\sim 2^{k_{\rm max}} \geq 2^{2j_{1}+j}$, this case can be proved
similarly to $2^{k_{2}}\sim 2^{k_{\rm max}}$ of (ii).

(c) When $2^{k}\geq 2^{2j_{1}+j}.$ In this case, $|\xi|\left\langle\tau-\xi^{3}+\frac{1}{\xi}
\right\rangle^{-1}\leq C 2^{-2j_{1}}$ or $|\xi|\leq 2^{k-2j_{1}}$ in $B_{k}.$

By using the H\"older inequality in $\xi$ and the Young inequality as well as (\ref{2.038}),
we have that
\begin{eqnarray}
&&\left\||\xi|\left\langle\tau-\xi^{3}+\frac{1}{\xi}\right\rangle^{-1}\langle\xi\rangle^{-3/4}f*g\right\|_{L_{\xi}^{2}
L_{\tau}^{1}(A_{0})}\nonumber\\&&\leq C2^{-2j_{1}}\|f*g\|_{L_{\xi}^{\infty}L_{\tau}^{1}
(A_{0})}\nonumber\\&&\leq C2^{-j_{1}/2}\left\|\langle\xi\rangle^{-3/4}f\right\|_{L_{\xi}^{2}L_{\tau}^{1}}
\left\|\langle\xi\rangle^{-3/4}g\right\|_{L_{\xi}^{2}L_{\tau}^{1}}\nonumber\\&&
\leq C2^{-j_{1}/2}\|f\|_{\hat{X}}
\|g\|_{\hat{X}}\label{3.018}.
\end{eqnarray}
When $(\xi,\tau)\in D$,  since $|\xi|\left\langle\tau-\xi^{3}+\frac{1}{\xi}\right\rangle^{-1/2}\leq C |\xi|^{1/2}2^{-j_{1}}$, by using Lemma 2.3, we have that
\begin{eqnarray}
&&\left\|I_{A_{j}}\left\langle\tau-\xi^{3}+\frac{1}{\xi}\right\rangle^{-1}\xi f\ast g\right\|_{\hat{X}^{-3/4,1/2}}\nonumber\\&&\leq
C\left\|\xi\left\langle\tau-\xi^{3}+\frac{1}{\xi}\right\rangle^{-1/2}f*g\right\|_{L^{2}}\nonumber
\\&&\leq C2^{-j_{1}}\left\||\xi|^{1/2}f*g\right\|_{L^{2}}\nonumber\\
&&\leq C2^{-j_{1}}2^{-\frac{j_{1}}{2}}\|f\|_{\hat{X}^{0,1/2,1}}
\|g\|_{\hat{X}^{0,1/2,1}}\nonumber\\&&\leq C\|f\|_{\hat{X}^{-3/4,1/2,1}}\|g\|_{\hat{X}^{-3/4,1/2,1}}\nonumber\\&&
\leq C\|f\|_{\hat{X}}
\|g\|_{\hat{X}}.\label{3.019}
\end{eqnarray}
When $(\xi,\tau)$ is  outside of  $D$ and $|\xi|\leq \frac{1}{8}$, we have that
\begin{eqnarray}
|\xi|^{-3}\geq |\tau|=\left|\tau-\xi^{3}+\frac{1}{\xi}\right|+1-|\xi|^{3}
-|\xi|^{-1}-1\geq 2^{k}-|\xi|^{3}-|\xi|^{-1}-1.\label{3.020}
\end{eqnarray}
From (\ref{3.020}), we have that $|\xi|\leq C2^{-j_{1}/2}.$
We consider the following two cases:
\begin{eqnarray}
&&(1): 2^{-(3-\epsilon)j_{1}}\leq |\xi|\leq C2^{-j_{1}/2},\label{3.021}\\
&&(2):|\xi|\leq 2^{-(3-\epsilon)j_{1}}\label{3.022}.
\end{eqnarray}
Case (1) can be proved similarly to $f*g\subset \left\{(\xi,\tau)\in \R^{2}:|\xi|\leq 1,
|\tau|\leq |\xi|^{-3}\right\}$ of (iv) in  \cite{Kis}.
Now we deal with case (2). In this case, we have that
\begin{eqnarray}
&&\left\|I_{A_{j}}\left\langle\tau-\xi^{3}+\frac{1}{\xi}\right
\rangle^{-1}\xi f*g\right\|_{\hat{X}}\nonumber\\&&\leq
C\left\|I_{A_{j}}\left\langle\tau-\xi^{3}+\frac{1}{\xi}\right
\rangle^{-1}\xi f*g\right\|_{\hat{X}^{-3/4,1/2,1}}\nonumber\\
&&\leq C2^{j}\sum_{k\geq0}2^{-k/2}\|f*g\|_{L^{2}(A_{j}\cap B_{k})}
\leq C(v_{1}+v_{2}),\label{3.023}
\end{eqnarray}
where
\begin{eqnarray}
&&v_{1}=\sum_{j=-\infty}^{-(3-\epsilon)j_{1}}2^{j}\sum_{k\geq2j_{1}+j}^{10j_{1}}
2^{-k/2}\left\|f*g\right\|_{L^{2}(B_{k})},\nonumber\\
&&v_{2}=\sum_{j=-\infty}^{-(3-\epsilon)j_{1}}2^{j}\sum_{k\geq10j_{1}} 2^{-k/2}
\left\|f*g\right\|_{L^{2}(B_{k})}\nonumber.
\end{eqnarray}
By using the Cauchy-Schwartz inequality with respect to $\tau$ and $\frac{j_{1}}{2^{(1-\epsilon)j_{1}}}\leq 1(j\geq 10)$, we have that
\begin{eqnarray}
&&v_{1}\leq Cj_{1}\sum_{j=-\infty}^{-(3-\epsilon)j_{1}}2^{j}\|f*g\|_{L_{\xi}^{2}L_{\tau}^{\infty}}
\leq Cj_{1}\sum_{j=-\infty}^{-(3-\epsilon)j_{1}}2^{j}\|f\|_{L_{\xi}^{2}L_{\tau}^{2}}
\|g\|_{L_{\xi}^{1}L_{\tau}^{2}}\nonumber\\&&\leq C
j_{1}2^{-(\frac{5}{2}-\epsilon)j_{1}}\|f\|_{L_{\xi}^{2}L_{\tau}^{2}}\|g\|_{L_{\xi}^{2}
L_{\tau}^{2}}\leq Cj_{1}2^{-(1-\epsilon)j_{1}}\|f\|_{\hat{X}^{-3/4,1/2,1}}
\|g\|_{\hat{X}^{-3/4,1/2,1}}\nonumber\\&&\leq C\|f\|_{\hat{X}^{-3/4,1/2,1}}
\|g\|_{\hat{X}^{-3/4,1/2,1}}\nonumber\\&&\leq C\|f\|_{\hat{X}}
\|g\|_{\hat{X}}\nonumber.
\end{eqnarray}
By using Lemma 2.3, we have that
\begin{eqnarray}
&&v_{2}\leq C2^{-11j_{1}/2}\sum_{j=-\infty}^{-(3-\epsilon)j_{1}}2^{j/2}
\|f\|_{\hat{X}^{0,1/2,1}}
\|g\|_{\hat{X}^{0,1/2,1}}\leq C2^{-\frac{14-\epsilon}{2}j_{1}}
\|f\|_{\hat{X}^{0,1/2,1}}
\|g\|_{\hat{X}^{0,1/2,1}}\nonumber\\&&\leq C2^{-\frac{11-\epsilon}{2}j_{1}}
\|f\|_{\hat{X}^{-3/4,1/2,1}}
\|g\|_{\hat{X}^{-3/4,1/2,1}}\leq C\|f\|_{\hat{X}^{-3/4,1/2,1}}
\|g\|_{\hat{X}^{-3/4,1/2,1}}\nonumber\\&&\leq C\|f\|_{\hat{X}}
\|g\|_{\hat{X}}\nonumber.
\end{eqnarray}
When $(\xi,\tau)$ is  outside of  $D$ and $|\xi|\geq \frac{1}{8}$,
this case can be proved similarly to (ii).

\noindent(vii) When $(\xi_{2},\tau_{2})\in A_{0}$ we consider $|\xi_{2}|\leq 2^{-2j},$
 $2^{-2j}\leq|\xi_{2}|\leq \frac{1}{8} $ and $\frac{1}{8}<|\xi_{2}|\leq 1$, respectively.

\noindent (a) Case $|\xi_{2}|\leq 2^{-2j}.$
By using the Young inequality and the H\"older inequality and (\ref{2.072}), we have that
\begin{eqnarray}
&&\left\|I_{A_{j}}\xi f*g\right\|_{\hat{X}^{-3/4,-1/2,1}}\leq C2^{j}\left\|\langle\xi\rangle^{-3/4}f
\right\|_{L_{\xi}^{2}L_{\tau}^{1}}\left\|g\right\|_{L_{\xi}^{1}L_{\tau}^{2}}\nonumber\\&&\leq
C\left\|\langle\xi\rangle^{-3/4}f\right\|_{L_{\xi}^{2}L_{\tau}^{1}}
\left\|g\right\|_{L_{\xi}^{2}L_{\tau}^{2}}\nonumber\\&&\leq C\|f\|_{\hat{X}^{-3/4,1/2,1}}
\left\|g\right\|_{\hat{X}^{-3/4,1/2}}\nonumber\\&&\leq C\|f\|_{\hat{X}}
\|g\|_{\hat{X}}.\label{3.024}
\end{eqnarray}
(b) $2^{-2j}\leq |\xi_{2}|\leq\frac{1}{8} ,$
cases  $2^{k}\sim 2^{k_{max}}$ and $2^{k_{1}}\sim 2^{k_{max}}$ can be proved similarly to
cases  $2^{k}\sim 2^{k_{max}}$ and $2^{k_{1}}\sim 2^{k_{max}}$ of (v) in \cite{Kis}.

\noindent Case  $2^{k_{2}}\sim 2^{k_{max}}$ and $(\tau_{2},\xi_{2})\in D.$

\noindent We only consider $2^{k_{2}}\sim 2^{k_{max}}>4{\rm max}
\left\{2^{k},2^{k_{1}}\right\},$ otherwise,
$2^{k_{2}}\sim 2^{k_{max}}\leq 4{\rm max}\left\{2^{k},2^{k_{1}}\right\}$
which boils down to
cases  $2^{k}\sim 2^{k_{max}}$ and $2^{k_{1}}\sim 2^{k_{max}}$.

In this case, we claim that $3|\xi\xi_{1}\xi_{2}|\geq \frac{\xi_{1}^{2}+\xi_{1}\xi_{2}
+\xi_{2}^{2}}{|\xi\xi_{1}\xi_{2}|}$.
If  $3|\xi\xi_{1}\xi_{2}|\leq \frac{\xi_{1}^{2}+\xi_{1}\xi_{2}+\xi_{2}^{2}}{|\xi\xi_{1}\xi_{2}|},$
since $|\tau_{2}|\geq |\xi_{2}|^{-3},$ then we have
\begin{eqnarray}
\frac{1}{2}|\xi_{2}|^{-3}\leq |\tau_{2}|-|\xi_{2}|^{3}-\frac{1}{|\xi_{2}|}
\leq\left|\tau_{2}-\xi_{2}^{3}+\frac{1}{\xi_{2}}\right|\leq \frac{4(\xi_{1}^{2}
+\xi_{1}\xi_{2}+\xi_{2}^{2})}{|\xi\xi_{1}\xi_{2}|}\label{3.025},
\end{eqnarray}
from (\ref{3.025}), we have
\begin{eqnarray}
\frac{\xi_{1}^{2}+\xi_{1}\xi_{2}+\xi_{2}^{2}}{|\xi\xi_{1}|}\geq \frac{1}{8}|\xi_{2}|^{-2}\geq 8\label{3.026},
\end{eqnarray}
from (\ref{3.026}), since $\xi\xi_{1}>0,$ we have that
\begin{eqnarray}
\xi_{1}^{2}+\xi_{1}\xi_{2}+\xi_{2}^{2}\geq 8|\xi\xi_{1}|=8\xi\xi_{1}=8\xi_{1}^{2}
+8\xi_{1}\xi_{2},\label{3.027}
\end{eqnarray}
which yields
\begin{eqnarray}
7\xi_{1}^{2}+7\xi_{1}\xi_{2}\leq \xi_{2}^{2}\leq \frac{1}{64}\label{3.028}.
\end{eqnarray}
(\ref{3.028}) is invalid. Thus, $3|\xi\xi_{1}\xi_{2}|\geq \frac{\xi_{1}^{2}
+\xi_{1}\xi_{2}+\xi_{2}^{2}}{|\xi\xi_{1}\xi_{2}|}$.
Thus, in this case $2^{-3j_{2}}\sim\frac{1}{2}|\xi_{2}|^{-3}\leq |\tau_{2}-\xi_{2}^{3}
+\frac{1}{\xi_{2}}|\sim |\xi\xi_{1}\xi_{2}|\sim 2^{2j+j_{2}}$
which yields $|\xi_{2}|\sim 2^{j_{2}}\sim2^{-j/2}.$
Consequently, $C2^{3j/2}\leq |\tau_{2}|\sim \left|\tau_{2}-\xi_{2}^{3}+\frac{1}{\xi_{2}}
\right|\leq C2^{2j}.$
Without loss of generality, we can assume that $g$ is supported on $A_{0}\cap B_{[3j/2,2j]}$.

\noindent(1)  When $g$ is supported on $[B_{3j/2,\frac{3j}{2}+\gamma]}$ with $0\leq \gamma \leq \frac{j}{2}$,
 for any
$\gamma^{'}\geq0,$ by using the Young inequality and $|\xi|\leq C|\xi_{2}|^{-1/2},$ we have that
\begin{eqnarray}
&&\left\|I_{B_{\geq \gamma^{'}}}\xi f*g\right\|_{\widehat{X}^{-3/4,-\frac{1}{2},1}}
\nonumber\\&&\leq C\sum_{k\geq \gamma^{'}}2^{-k/2}\left\|\langle\xi\rangle^{-3/4}f\right\|_{L_{\xi}^{2}
L_{\tau}^{1}}\left\||\xi|^{-1/2}\left\langle\tau-\xi^{3}
+\frac{1}{\xi}\right\rangle ^{1/2}g\right\|_{L_{\xi}^{1}L_{\tau}^{2}(B_{[\frac{3j}{2},\frac{3j}{2}+\gamma]})}\nonumber\\
&&\leq C2^{-\gamma^{'}/2}\|f\|_{\widehat{X}^{-3/4,-\frac{1}{2},1}}\||\xi|^{-1/2}\|_{L_{\xi}^{2}
(\{C2^{-j/2}\leq |\xi|\leq C2^{-\frac{j}{2}+\gamma}\})}\|g\|_{\widehat{X}^{0,\frac{1}{2}}}\nonumber\\
&&\leq C(\langle\gamma\rangle)^{1/2}2^{-\gamma^{'}/2}\|f\|_{\widehat{X}^{-3/4,\frac{1}{2},1}}
\|g\|_{\widehat{X}^{-3/4,\frac{1}{2}}}\nonumber\\&&\leq C\|f\|_{\hat{X}}
\|g\|_{\hat{X}}\label{3.029}
\end{eqnarray}
(2) When $g$ is supported on $B_{[\frac{3j}{2}+\gamma^{'},\frac{3j}{2}]}$ with $0\leq
 \gamma^{'} \leq
 \frac{j}{2}$, for any $\gamma\geq0,$ by using the Young inequality and $|\xi|\leq
 C|\xi_{2}|^{-1/2},$ by using Lemma 2.6,  we have that
\begin{eqnarray}
&&\left\|I_{B_{\leq \gamma}}\xi f*g\right\|_{\widehat{X}^{-3/4,-\frac{1}{2},1}}\nonumber
\\&&\leq C\sum_{k\leq \gamma}2^{-k/2}\left\|(\langle\xi\rangle^{-3/4}f)*(|\xi|^{-1/2}\left\langle\tau-\xi^{3}
+\frac{1}{\xi}\right\rangle ^{1/2}g)\right\|_{L^{2}(B_{k})}\nonumber\\
&&\leq C\langle\gamma\rangle 2^{-j/2}\|f\|_{\widehat{X}^{-3/4,\frac{1}{2},1}}\left\||\xi|^{-1/2}
\left\langle\tau-\xi^{3}+\frac{1}{\xi}\right\rangle^{1/2}g\right\|_{L_{\tau\xi}^{2}(|\xi|\geq
C2^{-\frac{j}{2}+\gamma^{'}})}\nonumber\\
&&\leq C(\langle\gamma\rangle)2^{-\gamma^{'}/2}\|f\|_{\widehat{X}^{-3/4,\frac{1}{2},1}}
\|g\|_{\widehat{X}^{-3/4,\frac{1}{2}}}\nonumber\\&&\leq C\|f\|_{\hat{X}}
\|g\|_{\hat{X}}\label{3.030}
\end{eqnarray}
Let $\gamma_{0}=\frac{j}{2}(\geq 2^{N+2})$, $\gamma_{n+1}=2{\rm log}_{2}{\gamma_{n}},$
$6\leq \gamma_{N}<8.$
Firstly, we apply (1) with $\gamma=\gamma_{0}$ and $\gamma^{'}=\gamma_{1}$, then apply (2)
with $\gamma=\gamma_{1}, \gamma^{'}=\gamma_{2}$. Repeating this procedure,  at the end
applying (1) with
$\gamma=\gamma_{N-1}$ and $\gamma^{'}=0,$ combining (1) with (2), by using Lemma \ref{Lemma2.7},
we have that
\begin{eqnarray}
&&\left\|\xi f*g\right\|_{\widehat{X}^{-3/4,-\frac{1}{2},1}}\nonumber\\&&\leq C(1+\sum_{n=0}^{N-1} \frac{1}{\gamma_{n}^{1/2}})\|f\|_{\widehat{X}^{-3/4,\frac{1}{2},1}}\|g\|_{\widehat{X}^{-3/4,\frac{1}{2},1}}\nonumber\\
&&\leq C\|f\|_{\widehat{X}^{-3/4,\frac{1}{2},1}}\|g\|_{\widehat{X}^{-3/4,\frac{1}{2}}}\nonumber\\&&\leq C\|f\|_{\hat{X}}
\|g\|_{\hat{X}}.\label{3.031}
\end{eqnarray}
 Case  $2^{k_{2}}\sim 2^{k_{max}}$ and $(\xi_{2},\tau_{2})$ is outside of $D.$
In this case, from Lemma \ref{Lemma2.8}, we have that $2^{k_{2}}\sim 2^{k_{\rm max}}\geq C 2^{2j}|\xi_{2}|$
which yields that $|\xi_{2}|\leq C2^{k_{2}-2j}$.
By using the proof similar to (3.10) of \cite{Kis}, we have that
\begin{eqnarray}
\left\|\xi f*g\right\|_{\widehat{X}^{-3/4,-\frac{1}{2},1}}
\leq C\|f\|_{\widehat{X}^{-3/4,\frac{1}{2},1}}\|g\|_{\widehat{X}^{-3/4,\frac{1}{2},1}}\leq C\|f\|_{\hat{X}}
\|g\|_{\hat{X}}.\label{3.032}
\end{eqnarray}
By using the Cauchy-Schwartz inequality in $\tau$, we have that (\ref{3.02}) can be bounded by
\begin{eqnarray*}
\left\|\xi f*g\right\|_{\widehat{X}^{-3/4,-\frac{1}{2},1}}
\end{eqnarray*}
in this case, which can be proved similarly to (\ref{3.01}) in this case.

\noindent (c) $\frac{1}{8}< |\xi_{2}|<1$. This case can be proved similarly to case
(iii).

\noindent(viii) This case can be proved similarly to case (vii).

We have completed the proof of Lemma 3.1.

\begin{Lemma}\label{Lem3.2}
Let  $u,v\in X$, then
\begin{eqnarray}
&&\left\|\mathscr{F}^{-1}\left[\left\langle\tau-\xi^{3}+
\frac {1}{\xi}\right\rangle^{-1}\mathscr{F}\left[\partial_{x}(uv)\right]\right]
\right\|_{X}+
\left\|\mathscr{F}^{-1}\left[\left\langle\tau-\xi^{3}+
 \frac{1}{\xi}\right\rangle^{-1}\mathscr{F}\left[\partial_{x}(uv)\right]\right]
\right\|_{Y}\nonumber\\&&\leq C\|u\|_{X}\|v\|_{X}.\label{3.045}
\end{eqnarray}
\end{Lemma}
{\bf Proof.} To prove (\ref{3.045}), it suffices to prove that
\begin{eqnarray}
&&\left\|\mathscr{F}^{-1}\left[\left\langle\tau-\xi^{3}+
\frac {1}{\xi}\right\rangle^{-1}\mathscr{F}\left[\partial_{x}(uv)\right]\right]
\right\|_{X}\leq C\|u\|_{X}\|v\|_{X}.\label{3.046}\\
&&\left\|\mathscr{F}^{-1}\left[\left\langle\tau-\xi^{3}+
 \frac{1}{\xi}\right\rangle^{-1}\mathscr{F}\left[\partial_{x}(uv)\right]\right]
\right\|_{Y}\leq C\|u\|_{X}\|v\|_{X}.\label{3.047}
\end{eqnarray}
We first prove (\ref{3.046}). By using $\|f\|_{\hat{X}}^{2}=\sum\limits_{j\geq 0}\|I_{A_{j}}f\|_{\hat{X}}^{2}$,
we have that
\begin{eqnarray*}
&&\left\|\mathscr{F}^{-1}\left[\left\langle\tau-\xi^{3}+
\frac {1}{\xi}\right\rangle^{-1}\mathscr{F}\left[\partial_{x}(uv)\right]\right]
\right\|_{X}^{2}\nonumber\\&&=\sum_{j,j_{1}j_{2}\geq0}\left\|\xi\left\langle\tau-\xi^{3}+
\frac {1}{\xi}\right\rangle^{-1}I_{A_{j}}(I_{A_{j_{1}}}\mathscr{F}u)*(I_{A_{j_{2}}}\mathscr{F}v)
\right\|_{\hat{X}}^{2}=\sum_{j=1}^{6}T_{j},
\end{eqnarray*}
where
\begin{eqnarray*}
&&T_{1}=\sum_{j,j_{1}j_{2}\geq0,\>i}\left\|\xi\left\langle\tau-\xi^{3}+
\frac {1}{\xi}\right\rangle^{-1}I_{A_{j}}(I_{A_{j_{1}}}\mathscr{F}u)*(I_{A_{j_{2}}}\mathscr{F}v)
\right\|_{\hat{X}}^{2},\nonumber\\
&&T_{2}=\sum_{j,j_{1}j_{2}\geq0,\>ii}\left\|\xi\left\langle\tau-\xi^{3}+
\frac {1}{\xi}\right\rangle^{-1}I_{A_{j}}(I_{A_{j_{1}}}\mathscr{F}u)*(I_{A_{j_{2}}}\mathscr{F}v)
\right\|_{\hat{X}}^{2},\nonumber\\
&&T_{3}=\sum_{j,j_{1}j_{2}\geq0,\>iii}\left\|\xi\left\langle\tau-\xi^{3}+
\frac {1}{\xi}\right\rangle^{-1}I_{A_{j}}(I_{A_{j_{1}}}\mathscr{F}u)*(I_{A_{j_{2}}}\mathscr{F}v)
\right\|_{\hat{X}}^{2},\nonumber\\
&&T_{4}=\sum_{j,j_{1}j_{2}\geq0,\>iv}\left\|\xi\left\langle\tau-\xi^{3}+
\frac {1}{\xi}\right\rangle^{-1}I_{A_{j}}(I_{A_{j_{1}}}\mathscr{F}u)*(I_{A_{j_{2}}}\mathscr{F}v)
\right\|_{\hat{X}}^{2},\nonumber\\
&&T_{5}=\sum_{j,j_{1}j_{2}\geq0,\>v}\left\|\xi\left\langle\tau-\xi^{3}+
\frac {1}{\xi}\right\rangle^{-1}I_{A_{j}}(I_{A_{j_{1}}}\mathscr{F}u)*(I_{A_{j_{2}}}\mathscr{F}v)
\right\|_{\hat{X}}^{2},\nonumber\\
&&T_{6}=\sum_{j,j_{1}j_{2}\geq0,\>vi}\left\|\xi\left\langle\tau-\xi^{3}+
\frac {1}{\xi}\right\rangle^{-1}I_{A_{j}}(I_{A_{j_{1}}}\mathscr{F}u)*(I_{A_{j_{2}}}\mathscr{F}v)
\right\|_{\hat{X}}^{2},\nonumber\\
&&T_{7}=\sum_{j,j_{1}j_{2}\geq0,\>vii}\left\|\xi\left\langle\tau-\xi^{3}+
\frac {1}{\xi}\right\rangle^{-1}I_{A_{j}}(I_{A_{j_{1}}}\mathscr{F}u)*(I_{A_{j_{2}}}\mathscr{F}v)
\right\|_{\hat{X}}^{2},\nonumber\\
&&T_{8}=\sum_{j,j_{1}j_{2}\geq0,\>viii}\left\|\xi\left\langle\tau-\xi^{3}+
\frac {1}{\xi}\right\rangle^{-1}I_{A_{j}}(I_{A_{j_{1}}}\mathscr{F}u)*(I_{A_{j_{2}}}\mathscr{F}v)
\right\|_{\hat{X}}^{2}.
\end{eqnarray*}
here $(i),(ii),(iii),(iv),(v),(vi),(vii), (viii)$ is  case $(i),(ii),(iii),(iv),(v),(vi),(vii), (viii)$ of Lemmas 3.1, 3.2.
Combining $T_{j} (1\leq j\leq 8,j\in N)$, Lemmas 3.1, 3.2 with $\|f\|_{\hat{X}}^{2}=\sum_{j\geq 0}\|I_{A_{j}}f\|_{\hat{X}}^{2}$, we easily obtain
(\ref{3.046}). By using a proof similarly to (\ref{3.046}), we easily obtain (\ref{3.047}).

We have completed the proof of Lemma 3.2.

\bigskip
\bigskip

\noindent {\large\bf 4. Proof of Theorem  1.1}

\setcounter{equation}{0}

 \setcounter{Theorem}{0}

\setcounter{Lemma}{0}

\setcounter{section}{4}

(\ref{1.01})-(\ref{1.02}) is equivalent to the following integral equation:
\begin{eqnarray}
u(t)=e^{-t(-\partial_{x}^{3}+\partial_{x}^{-1})}u_{0}+\frac{1}{2}\int_{0}^{t}e^{-(t-s)
(-\partial_{x}^{3}+\partial_{x}^{-1})}\partial_{x}(u^{2})ds.\label{5.01}
\end{eqnarray}
We define
\begin{eqnarray}
\Phi(u)=e^{-t(-\partial_{x}^{3}+\partial_{x}^{-1})}u_{0}+\frac{1}{2}\int_{0}^{t}e^{-(t-s)
(-\partial_{x}^{3}+\partial_{x}^{-1})}\partial_{x}(u^{2})ds.\label{5.02}
\end{eqnarray}
By using Lemma \ref{Lemma2.10} and Lemma \ref{Lem3.2}, we have that
\begin{eqnarray}
\|\Phi(u)\|_{X_{1}}+\sup\limits_{-1\leq t\leq1}\|\Phi(u)\|_{H^{-3/4}(\SR)}
\leq C\|u_{0}\|_{H^{-3/4}(\SR)}+C\|u\|_{X_{1}}^{2},
\end{eqnarray}
when $\|u_{0}\|_{H^{-3/4}}$ is sufficiently small, we have that $\Phi(u)$ is
a contraction mapping on some closed ball in $X_{1}\cap C_{t}^{0}([-1,1]; H^{-3/4}(\R))$.
Thus $\Phi$ have a fixed point $u$, which is the local solution of (\ref{5.01})
and thus (\ref{1.01})(\ref{1.02}). For large data, by taking
$u_{\lambda 0}(x)=\lambda^{-2}u_{0}\left(\frac{x}{\lambda}\right),$
we have that
$
  \|u_{\lambda 0}\|_{H^{-3/4}(\SR)}\leq C\lambda^{-3/4}\|u_{0}\|_{H^{-3/4}(\SR)}.
$
 Taking $\lambda$ sufficiently large, then
$\|u_{\lambda 0}\|_{H^{s}(\SR)}$  is sufficiently small, then
there is a  solution to (\ref{1.01}) associated to the initial function $u_{\lambda 0}(x,0)$,
and thus (\ref{1.01})(\ref{1.02}) admit a  solution.  The Lipschitz dependence of solutions
on the data and the uniqueness of the solutions can be found in  \cite{MT,Kis}.

We have completed the proof of Theorem 1.1.

\bigskip
\bigskip
\bigskip

\section{\large\bf  Appendix}

\setcounter{equation}{0}

 \setcounter{Theorem}{0}

\setcounter{Lemma}{0}

\setcounter{section}{5}

{\bf Example 1.} (high $\times$ high$\mapsto$ low interaction.)
Let $Rec$  be the region in $R^{2}_{\tau\xi}$ inside the parallelogram with vertices
\begin{eqnarray}
&&(\tau,\xi)=(N^{3},N), (N^{3}+N^{\frac{3}{2}},N+\frac{1}{3}N^{-\frac{1}{2}}),\label{7.01}\\
&&\left((N+\frac{1}{3}N^{-\frac{1}{2}})^{3}, N+\frac{1}{3}N^{-\frac{1}{2}}\right), \left(N^{3}+\frac{1}{3}+\frac{1}{27}N^{-\frac{3}{2}},N\right),\label{6.02}
\end{eqnarray}
where $N$ is a sufficiently large positive number. It is easily checked that
$Rec$ is included in the region $\left\{|\tau-\xi^{3}+\frac{1}{\xi}|<1\right\},$ has
the longest side pointing at the direction $(3N^{2},1)$ and $|Rec|\sim N^{-1/2}.$ We put
$R_{0}$ equal to the translation of $R$ centered at the origin.
Let
\begin{eqnarray}
\mathscr{F}u(\tau,\xi):=I_{Rec}, \mathscr{F}v(\tau,\xi)=\mathscr{F}u(-\tau,-\xi),
\end{eqnarray}
where $I_{\Omega}$ denotes the characteristic function of a set $\Omega.$ By a direct
computation, we have that
\begin{eqnarray}
\|u\|_{X^{-\frac{3}{4},b}}=\|v\|_{X^{-\frac{3}{4},b}}\sim N^{-1}, \mathscr{F}(uv)\geq CN^{-\frac{1}{2}}I_{R_{0}},
\|\partial_{x}(uv)\|_{X^{-\frac{3}{4},b-1}}\geq CN^{\frac{6b-11}{4}}.\label{6.03}
\end{eqnarray}
Then
\begin{eqnarray}
\|\partial_{x}(fg)\|_{X^{-\frac{3}{4},\>b-1}}\leq C\|f\|_{X^{-\frac{3}{4},\>b}}\|g\|_{X^{-\frac{3}{4},\>b}}\label{6.04}.
\end{eqnarray}
is invalid for $b>1/2.$

{\bf Example 2.} (high $\times$ low $\mapsto$ high interaction.) Let $\mathscr{F}u=I_{Rec},$
$\mathscr{F}v=I_{R_{0}}.$ By a direct computation, we have that
\begin{eqnarray}
&&\|u\|_{X^{-\frac{3}{4},b}}\sim N^{-1},\|v\|_{X^{-\frac{3}{4},b}}\sim N^{\frac{6b-1}{4}}, \nonumber\\&&\mathscr{F}(uv)
\geq CN^{-\frac{1}{2}}I_{Rec},\|\partial_{x}(uv)\|_{X^{-frac{3}{4},b}}\geq CN^{1/2},\label{6.05}
\end{eqnarray}
thus,
\begin{eqnarray}
\|\partial_{x}(fg)\|_{X^{-\frac{3}{4},\>b-1}}\leq C\|f\|_{X^{-\frac{3}{4},\>b}}\|g\|_{X^{-\frac{3}{4},\>b}}\label{6.06}
\end{eqnarray}
is invalid for $b<1/2.$

\bigskip

\leftline{\large \bf Acknowledgments}

\bigskip

\noindent

 This work is supported by the Natural Science Foundation of China
 under grant numbers 11171116 and 11401180. The first author is also
 supported in part by the Fundamental Research Funds for the
 Central Universities of China under the grant number 2012ZZ0072.
  The second author is  supported by the
 NSF of China (No.11371367) and Fundamental research program of
  NUDT(JC12-02-03).
 The third author is also supported by the Natural Science
 Foundation of China
 under grant number 14IRTSTHN023..

\end{document}